\documentclass[12pt,amscd]{amsart}
\usepackage{amssymb}
\numberwithin{equation}{section}
\newtheorem{thm}[equation]{Theorem}

\newtheorem{cor}[equation]{Corollary}

\newtheorem{lem}[equation]{Lemma}

\newtheorem{prop}[equation]{Proposition}

\newenvironment{pf}{\proof[\proofname]}{\endproof}
\newenvironment{pf*}[1]{\proof[#1]}{\endproof}

\theoremstyle{definition}

\newtheorem{defn}[equation]{Definition}

\theoremstyle{remark}

\newtheorem*{rmk}{Remark}
\newtheorem*{rmks}{Remarks}

\newtheorem*{ack}{Acknowledgement}

\newcommand{\comment}[1]{}

\makeatother

\begin{document}
\baselineskip=18truept


\def\C {{\mathbb C}}
\def\Cn {{\mathbb C}^n}
\def\R {{\mathbb R}}
\def\Rn {{\mathbb R}^n}
\def\Z {{\mathbb Z}}
\def\N {{\mathbb N}}
\def\cal#1{{\mathcal #1}}
\def\bb#1{{\mathbb #1}}

\def\dbar {\bar \partial }
\def\dir {{\mathcal D}}
\def\lev#1{{\mathcal L}\left(#1\right)}
\def\lap {\Delta }
\def\ol {{\mathcal O}}
\def\E {{\mathcal E}}
\def\J {{\mathcal J}}
\def\U {{\mathcal U}}
\def\V {{\mathcal V}}
\def\z {\zeta }
\def\Harm {\text {Harm}\, }
\def\grad {\nabla }
\def\dexh {\{ M_k \} _{k=0}^{\infty } }
\def\sing#1{#1_{\text {sing}}}
\def\reg#1{#1_{\text {reg}}}
\def\pione#1{\pi_1(#1)}
\def\pioneat#1#2{\pi_1(#1,#2)}

\def\lquotient{\big{\backslash}}

\def\setof#1#2{\{ \, #1 \mid #2 \, \} }

\def\image#1{\text{\rm im}\,\bigl[#1\bigr]}
\def\kernel#1{\text{\rm ker}\,\bigl[#1\bigr]}

\def\holecl {M\setminus \overline M_0}
\def\hole {M\setminus M_0}

\def\nd{\frac {\partial }{\partial\nu } }
\def\ndof#1{\frac {\partial#1}{\partial\nu } }

\def\pdof#1#2{\frac {\partial#1}{\partial#2}}
\def\pdwrt#1{\frac{\partial}{\partial #1}}

\def\cinf{\ensuremath{\mathcal C^{\infty}} }
\def\cinfns{\ensuremath{\mathcal C^{\infty}}}

\def\dist{\text{\rm dist}}

\def\diam{\text{\rm diam}}

\def\real{\text{\rm Re}\, }

\def\imag{\text{\rm Im}\, }

\def\supp{\text{\rm supp}\, }

\def\Vol{\text{\rm vol}}

\def\restrict#1{{\upharpoonright_{#1}}}

\def\sm{\setminus}

\def\plshclass{{\mathcal {P}}}
\def\strplshclass{{\mathcal S\mathcal P}}

\def\geqtrace#1#2{{\geq_{(#1,#2)}}}
\def\gtrace#1#2{{>_{(#1,#2)}}}
\def\leqtrace#1#2{{\leq_{(#1,#2)}}}
\def\ltrace#1#2{{<_{(#1,#2)}}}



\def\anal{analytic }
\def\analns{analytic}

\def\bdd{bounded }
\def\bddns{bounded}

\def\cpt{compact }
\def\cptns{compact}

\def\cpx{complex }
\def\cpxns{complex}

\def\cont{continuous }
\def\contns{continuous}

\def\dime{dimension }
\def\dimens{dimension }

\def\exh{exhaustion }
\def\exhns{exhaustion}

\def\fn{function }
\def\fnns{function}

\def\fns{functions }
\def\fnsns{functions}

\def\holo{holomorphic }
\def\holons{holomorphic}

\def\mero{meromorphic }
\def\merons{meromorphic}

\def\holoconvex{holomorphically convex }
\def\holoconvexns{holomorphically convex}

\def\ircomp{irreducible component }
\def\concomp{connected component }
\def\ircompns{irreducible component}
\def\concompns{connected component}
\def\ircomps{irreducible components }
\def\concomps{connected components }
\def\ircompsns{irreducible components}
\def\concompsns{connected components}

\def\irred{irreducible }
\def\irredns{irreducible}

\def\con{connected }
\def\conns{connected}

\def\comp{component }
\def\compns{component}
\def\comps{components }
\def\compsns{components}

\def\mfld{manifold }
\def\mfldns{manifold}
\def\mflds{manifolds }
\def\mfldsns{manifolds}

\def\nbd{neighborhood }
\def\nbds{neighborhoods }
\def\nbdns{neighborhood}
\def\nbdsns{neighborhoods}

\def\harm{harmonic }
\def\harmns{harmonic}
\def\plh{pluriharmonic }
\def\plhns{pluriharmonic}
\def\plsh{plurisubharmonic }
\def\plshns{plurisubharmonic}

\def\qplsh#1{$#1$-plurisubharmonic}
\def\hplsh{$(n-1)$-plurisubharmonic }
\def\hplshns{$(n-1)$-plurisubharmonic}

\def\para{parabolic }
\def\parans{parabolic}

\def\rel{relatively }
\def\relns{relatively}

\def\str{strictly }
\def\strns{strictly}

\def\strg{strongly }
\def\strgns{strongly}

\def\cvx{convex }
\def\cvxns{convex}

\def\wrt{with respect to }
\def\wrtns{with respect to}

\def\st {such that }
\def\stns {such that}

\def\hm {harmonic measure }
\def\hmns {harmonic measure}

\def\hmib {harmonic measure of the ideal boundary of }
\def\hmibns {harmonic measure of the ideal boundary of}

\def\Vert{\text{{\rm Vert}}\, }
\def\Edge{\text{{\rm Edge}}\, }

\def\til#1{\tilde{#1}}
\def\wtil#1{\widetilde{#1}}

\def\what#1{\widehat{#1}}

\def\seq#1#2{\{#1_{#2}\} }


\def\vphi {\varphi }


\def\inv{   ^{-1}  }

\def\ssp#1{^{(#1)}}

\def\set#1{\{ #1 \}}

\title[The Bochner--Hartogs dichotomy]
{The Bochner--Hartogs dichotomy for bounded geometry hyperbolic
K\"ahler manifolds}
\author[T.~Napier]{Terrence Napier}
\address{Department of Mathematics\\Lehigh University\\Bethlehem, PA 18015\\USA}
\email{tjn2@lehigh.edu}
\thanks{To appear in Annales de l'Institut Fourier.}
\author[M.~Ramachandran]{Mohan Ramachandran}
\address{Department of Mathematics\\University at Buffalo\\Buffalo, NY 14260\\USA}
\email{ramac-m@buffalo.edu}


\date{June 12, 2015}

\begin{abstract}
The main result is that for a \con hyperbolic complete K\"ahler
manifold with bounded geometry of order two and exactly one end,
either the first compactly supported cohomology with values in the
structure sheaf vanishes or the manifold admits a proper \holo
mapping onto a Riemann surface. 
\end{abstract}

\maketitle

\section*{Introduction} \label{introduction}

Let $(X,g)$ be a \con non\cpt complete K\"ahler manifold.
According to \cite{Gro-Sur la groupe fond}, \cite{Li Structure
complete Kahler}, \cite{Gro-Kahler hyperbolicity},
\cite{Gromov-Schoen}, \cite{NR-Structure theorems},
\cite{Delzant-Gromov Cuts}, \cite{NR Filtered ends}, and \cite{NR L2 Castelnuovo}, if $X$
has at least three filtered ends relative to the universal
covering (i.e., $\tilde e(X)\geq 3$ in the sense of Definition~\ref{ends filtered ends def}) and $X$ is weakly
$1$-complete (i.e., $X$ admits a \cont \plsh \exh \fnns) or $X$ is regular hyperbolic (i.e., $X$ admits a
positive symmetric Green's \fn that vanishes at infinity) or $X$
has \bdd geometry of order two (see Definition~\ref{Bdd geom along set defn}), then $X$ admits a proper \holo mapping onto a
Riemann surface. In particular, if $X$ has at least three
(standard) ends (i.e., $e(X)\geq 3$) and $X$ satisfies one of the above three
conditions, then such a mapping exists.
Cousin's example~\cite{Cousin} of a $2$-ended weakly $1$-complete
covering of an Abelian variety that has only constant \holo \fns
demonstrates that two (filtered) ends do not suffice.

A non\cpt \cpx manifold~$X$ for which $H^1_c(X,\ol)=0$ is said to
have the {\it Bochner--Hartogs property} (see
Hartogs~\cite{Hartogs}, Bochner~\cite{Bochner}, and
Harvey~and~Lawson~\cite{Harvey-Lawson}). Equivalently, for every
\cinf compactly supported form~$\alpha$ of type $(0,1)$ with
$\dbar\alpha=0$ on~$X$, there is a \cinf compactly supported \fn
$\beta$ on~$X$ \st $\dbar\beta=\alpha$. If $X$ has the
Bochner--Hartogs property, then every \holo \fn on a \nbd of
infinity with no \rel \cpt \concomps extends to a \holo \fn
on~$X$. For cutting off away from infinity, one gets a \cinf \fn
$\lambda$ on $X$. Taking $\alpha\equiv\dbar\lambda$ and forming
$\beta$ as above, one then gets the desired
extension $\lambda-\beta$. In particular, $e(X)=1$, since for a \cpx manifold with
multiple ends, there exists a locally constant \fn on a \nbd of
$\infty$ that is equal to~$1$ along one end and $0$ along the
other ends, and such a \fn cannot extend holomorphically. Thus in
a sense, the space $H^1_c(X,\ol)$ is a function-theoretic
approximation of the set of (topological) ends of~$X$.
An open Riemann surface~$S$, as well as any \cpx manifold
admitting a proper \holo mapping onto~$S$, cannot have the
Bochner--Hartogs property, because $S$ admits
meromorphic \fns with finitely many poles. Examples of manifolds
of dimension~$n$ having the Bochner--Hartogs property include
strongly $(n-1)$-complete \cpx manifolds (Andreotti and
Vesentini~\cite{Andreotti-Vesentini}) and strongly
hyper-$(n-1)$-convex K\"ahler manifolds (Grauert and
Riemenschneider~\cite{Grauert-Riemenschneider}). We will say that
the \emph{Bochner--Hartogs dichotomy} holds for a class of \con
\cpx manifolds if each element either has the Bochner--Hartogs
property or admits a proper \holo mapping onto a Riemann surface.

According
to \cite{Ramachandran-BH for coverings}, \cite{NR-BH Weakly
1-complete}, and \cite{NR-BH Regular hyperbolic Kahler}, the Bochner--Hartogs dichotomy
holds for the class of
weakly $1$-complete or regular hyperbolic complete K\"ahler manifolds with
exactly one end. The main goal of this paper is the following:
\begin{thm}\label{BHD bounded geom hyperbolic thm}
Let $X$ be a \con non\cpt hyperbolic complete K\"ahler manifold
with \bdd geometry of order two, and assume that $X$ has exactly
one end. Then $X$ admits a proper \holo mapping onto a Riemann
surface if and only if $H^1_c(X,\ol)\neq 0$.
\end{thm}
In other words, the Bochner--Hartogs dichotomy holds for the class
of hyperbolic \con non\cpt complete K\"ahler manifolds with \bdd
geometry of order two and exactly one end. When combined with the
earlier results, the above gives the following:
\begin{cor}\label{combined BHD cor}
Let $X$ be a \con non\cpt complete K\"ahler manifold that has exactly one end (or has
at least three filtered ends) and satisfies at least one of the following:
\begin{enumerate}

\item[(i)] $X$ is weakly $1$-complete;

\item[(ii)] $X$ is regular hyperbolic; or

\item [(iii)] $X$ is hyperbolic and of \bdd geometry of
order two.

\end{enumerate}
Then $X$ admits a proper \holo mapping onto a Riemann surface if
and only if $H^1_c(X,\ol)\neq 0$.
\end{cor}
In particular, since \con coverings of \cpt K\"ahler manifolds have bounded geometry of all orders, 
we have the following 
(cf.~\cite{ArapuraBressRam}, 
\cite{Ramachandran-BH for coverings}, and Theorem~0.2 of \cite{NR-BH Regular hyperbolic Kahler}):
\begin{cor}\label{BHD covering cor}
Let $X$ be a \cpt K\"ahler manifold, and $\what X\to X$ a \con infinite covering that is
hyperbolic and has exactly one end (or at least three filtered ends). Then $\what X$ admits a proper \holo mapping onto a Riemann
surface if and only if $H^1_c(\what X,\ol)\neq 0$.
\end{cor}

A standard method for constructing a proper \holo mapping onto
a Riemann surface is to produce suitable linearly independent
\holo $1$-forms (usually as \holo differentials of \plh \fnsns),
and to then apply versions of Gromov's cup product lemma and the
Castelnuovo--de~Franchis theorem. In this context, an
\emph{irregular} hyperbolic manifold has a surprising advantage
over a \emph{regular} hyperbolic manifold in that an irregular
hyperbolic complete K\"ahler manifold with bounded geometry of order two
automatically admits a nonconstant 
positive \plh \fnns.
In particular, the proof of Theorem~\ref{BHD bounded geom
hyperbolic thm} in the irregular hyperbolic case is, in a sense,
simpler than the proof in the regular hyperbolic case (which
appeared in~\cite{NR-BH Regular hyperbolic Kahler}).  Because
the existence of irregular
hyperbolic complete K\"ahler manifolds with one end and bounded
geometry of order two is not completely obvious, 
a $1$-dimensional example is provided in 
Section~\ref{BG irreg hyp example sect}.  However,  the authors do not know whether or not there exist examples with the above properties that
satisfy the Bochner--Hartogs property (and hence do not admit
proper \holo mappings onto Riemann surfaces).

Section~\ref{ends and bh basic sect} is a consideration of some
elementary properties of ends, as well as some
elementary topological properties of \cpx manifolds with the
Bochner--Hartogs property. Section~\ref{Bounded geometry sect}
contains the definition of bounded geometry. Section~\ref{greens
fn sect} consists of some terminology and facts from potential
theory, and a proof that the Bochner--Hartogs property holds for any
one-ended \con non\cpt hyperbolic complete K\"ahler manifold with
no nontrivial $L^2$ \holo $1$-forms. A modification of Nakai's
construction of an infinite-energy positive quasi-Dirichlet finite
\harm \fn on an irregular hyperbolic manifold, as well as a
modification of a theorem of Sullivan which gives pluriharmonicity
in the setting of a complete K\"ahler manifold with bounded
geometry of order two, appear in Section~\ref{quasidirichlet
finite fns sect}. The proof of Theorem~\ref{BHD bounded geom
hyperbolic thm} and the proofs of some related results appear in
Section~\ref{proof main thm sect}. An example of an irregular
hyperbolic complete K\"ahler manifold with one end and bounded
geometry of all orders is constructed in Section~\ref{BG irreg hyp
example sect}. 

 \begin{ack}
 The authors would like to thank the referee for providing 
 very valuable comments.
 \end{ack}

\section{Ends and the Bochner--Hartogs property}\label{ends and bh basic sect}

In this section, we consider an elementary topological property of
\cpx manifolds with the Bochner--Hartogs property. Further
topological characterizations of the Bochner--Hartogs dichotomy
will be considered in Section~\ref{proof main thm sect}. We first
recall some terminology and facts concerning ends.

By an {\it end} of a \con manifold~$M$, we will mean either a \comp $E$
of $M\setminus K$ with non\cpt closure, where $K$ is a given \cpt
subset of $M$, or an element of
\[
\lim _{\leftarrow} \pi _0 (M\setminus K),
\]
where the limit is taken as $K$ ranges over the \cpt subsets of
$M$ (or equivalently, the \cpt subsets of $M$ for which the complement $M\setminus K$
has no \rel \cpt \compsns, since the union of any \cpt subset of~$M$ with the \rel \cpt \concomps
of its complement is \cptns). The number of ends of $M$ will be
denoted by~$e(M)$. For a \cpt set $K$ \st $M\setminus K$ has no
\rel \cpt \compsns, we will call
\[
M\setminus K=E_1\cup \cdots \cup E_m,
\]
where $\seq Ej_{j=1}^m$ are the distinct \comps of $M\setminus K$,
an {\it ends decomposition} for~$M$.

\begin{lem}\label{basic ends lem}
Let $M$ be a \con non\cpt \cinf manifold.
\begin{enumerate}

\item[(a)] If $S\Subset M$,
then the number of \comps of $M\sm S$ that are not \rel \cpt in~$M$ is at most the number of
\comps of $M\sm T$ for
any set $T$ with $S\subset T\Subset M$. In
particular, the number of such \comps of $M\sm S$ is at most $e(M)$.

\item[(b)] If $K$ is a \cpt subset of~$M$, then there exists a
\cinf \rel \cpt domain~$\Omega$ in~$M$ containing~$K$ \st
$M\sm\Omega$ has no \cpt \compsns. In particular, if $k$ is a
positive integer with $k\leq e(M)$, then we may choose $\Omega$ so
that $M\sm\Omega$ also has at least $k$ \compsns; and hence
$\partial\Omega$ has at least $k$ \compsns.

\item[(c)] If $\Omega$ is a nonempty \rel \cpt domain in $M$, then the number of
\comps of~$\partial\Omega$ is at most $e(\Omega)$, with equality if $\Omega$ is also
smooth.

\item[(d)] Given an ends decomposition $M\sm K=E_1\cup\cdots\cup
E_m$, there is a \con \cpt set~$K'\supset K$ \st any
domain~$\Theta$ in~$M$ containing~$K'$ has an ends decomposition
$\Theta\sm K=E_1'\cup\cdots\cup E_m'$, where $E_j'=E_j\cap\Theta$
for $j=1,\dots,m$.

\item[(e)] If $\Omega$ and $\Theta$ are domains in~$M$ with
$\Theta\subset\Omega$ and both $M\sm\Omega$ and $\Omega\sm\Theta$
have no \cpt \compsns, then $M\sm\Theta$ has no \cpt \compsns.

\item[(f)] If $M$ admits a proper surjective \cont open mapping
onto an orientable topological surface that is not simply \conns,
then there exists a \cinf \rel \cpt domain~$\Omega$ in~$M$ \st
$M\sm\Omega$ has no \cpt \comps and $\partial\Omega$ is not
\conns.

\end{enumerate}
\end{lem}
\begin{pf}
For the proof of (a), we simply observe that if $S\subset T\Subset M$,
then each \concomp of $M\sm S$ that is not \rel \cpt in~$M$ must meet $M\sm T$ and must therefore contain some
\comp of $M\sm T$. Choosing $T\supset S$ to be a \cpt set for which $M\sm T$ has
no \rel \cpt \compsns, we see that the number of \comps of $M\sm S$ is at most
$e(M)$.

For the proof of (b), observe that given a \cpt set $K\subset M$, we may fix a \cinf
domain~$\Omega_0$  with $K\subset\Omega_0\Subset M$. The union of $\Omega_0$ with those (finitely many)
\comps of $M\sm\Omega_0$ which are \cpt is then a \cinf \rel \cpt domain~$\Omega\supset K$
in~$M$ for which $M\sm\Omega$ has no \cpt \compsns. Given a positive integer $k\leq e(M)$, we may choose~$\Omega$
to also contain a \cpt set~$K'$ for which $M\sm K'$ has at least~$k$ \comps and no \rel \cpt \comps in~$M$.
Part~(a) then implies that $M\sm\Omega$ has at least $k$ \compsns, and since each \comp must
contain a \comp of $\partial\Omega$, we see also that $\partial\Omega$ has at least~$k$ \compsns.

For the proof of (c), suppose $\Omega$ is a nonempty \rel \cpt domain in~$M$ and $k$ is a positive
integer. If $\partial\Omega$ has at least $k$ \compsns, then we may fix a covering of $\partial\Omega$
by disjoint \rel \cpt open
subsets $U_1,\dots,U_k$ of $M$ each of which meets $\partial\Omega$
(one may prove the existence of such sets by induction on~$k$).
We may also fix a \cpt set $K\subset\Omega$ containing
$\Omega\sm(U_1\cup\cdots\cup U_k)$ \st the \comps of $\Omega\sm
K\subset U_1\cup\cdots\cup U_k$ are not \rel \cpt in~$\Omega$. For
each~$j=1,\dots,k$, $U_j$ meets $\partial\Omega$ and therefore
some \compns~$E$ of $\Omega\sm K$, and hence $E\subset U_j$. Since
$\Omega\sm K$ has at most $e(\Omega)$ \compsns, it follows that
$k\leq e(\Omega)$. Furthermore, if $\Omega$ is smooth, then we may
choose $k$ to be equal to the number of boundary \compsns, and we
may choose the arbitrarily small \nbds so that $U_j\cap\Omega$ is
\con for each~$j$. We then get $k=e(\Omega)$ in this case.

For the proof of~(d), let $M\sm K=E_1\cup\cdots\cup E_m$ be an
ends decomposition. We may fix a \cinf \rel \cpt domain~$\Omega$
in~$M$ containing~$K$ \st $M\sm\Omega$ has no \cpt \compsns, and
for each $j=1,\dots,m$, we may fix a \con \cpt set~$A_j\subset
E_j$ \st $A_j$ meets each of the finitely many \comps of
$E_j\cap\partial\Omega$. The \cpt set
$K'\equiv\overline\Omega\cup\bigcup_{j=1}^mA_j$ then has the
required properties.

For the proof of (e), suppose $\Omega$ and $\Theta$ are domains
in~$M$ with $\Theta\subset\Omega$, and both $M\sm\Omega$ and
$\Omega\sm\Theta$ have no \cpt \compsns. If $E$ is a \comp of
$M\sm\Theta$, then either $E$ meets $M\sm\Omega$, in which case
$E$ contains a non\cpt \comp of $M\sm\Omega$, or $E\subset\Omega$,
in which case $E$ is a \comp of $\Omega\sm\Theta$.  In either
case, $E$ is non\cptns.

Finally for the proof of~(f), suppose $\Phi\colon M\to S$ is a
proper surjective \cont open mapping onto an orientable
topological surface~$S$ that is not simply \conns. By part~(b), we
may assume without loss of generality that $e(M)=1$. If $U$ is any
open set in~$S$ and $V$ is any \comp of $\Phi\inv(U)$, then
$\Phi(V)$ is both open and closed in~$U$; i.e., $\Phi(V)$ is a
\comp of~$U$. Consequently, if $K\subset S$ is a \cpt set for
which $S\sm K$ has no \rel \cpt \comps and $V$ is any \comp of
$M\sm\Phi\inv(K)=\Phi\inv(S\sm K)$, then $\Phi(V)$ must be a \comp
of $S\sm K$. Hence $V$ must be the unique \comp of
$M\sm\Phi\inv(K)$ that is not \rel \cpt in~$M$, and it follows
that $V=M\sm\Phi\inv(K)$ and $\Phi(V)=S\sm K$ are \conns. In
particular, $e(S)=1$.

Since every planar domain with one end is simply \conns, $S$ must
be nonplanar; that is, there exists a nonseparating simple closed
curve in~$S$.  Hence there exists a homeomorphism $\Psi$ of a suitable annulus
$\Delta(0;r',R')\equiv\setof{z\in\C}{r'<|z|<R'}$ onto a domain
$A'\Subset S$ with \con complement $S\sm A'$.  Fixing $r$ and $R$
with $0<r'<r<R<R'$, setting $A\equiv\Psi(\Delta(0;r,R))\Subset
A'$, $F\equiv S\sm\overline A$ and
$E\equiv\Phi\inv(F)=M\sm\Phi\inv(\overline A)$, and letting
$\Theta$ be a \comp of $\Phi\inv(A)$, we see that $E$ is \con and
$\Phi(\Theta)=A$ (by the above), and that $\overline E\subset
M\sm\Theta$. It also follows that $M\sm\Theta$ is \conns. For if
$K$ is a \cpt \comp of $M\sm\Theta$, then we must have
$K\cap\overline E=\emptyset$. Forming a \con \nbdns~$U$ of $K$ in
$M\sm\overline E\subset M\sm E=\Phi\inv(\overline A)$, we get
$\Phi(K)\subset\Phi(U)\subset\overline A$, and hence
$\Phi(U)\subset A$. Thus $K$ must lie in some
\compns~$V\subset M\sm\Theta$ of $\Phi\inv(A)$, and hence $K=V$.
But then $K$ must be both open and closed
in~$M$, which is clearly impossible. Therefore $M\sm\Theta$ is
\conns. Moreover, since $\Phi(\Theta)=A$, we must have
$\Phi(\partial\Theta)=\partial A$, and hence $\partial\Theta$ is
not \conns. Applying parts (b), (c), and~(e), we get the desired
smooth domain $\Omega\Subset\Theta$.
\end{pf}

As indicated in the introduction, a \con non\cpt \cpx manifold
with the Bochner--Hartogs property must have exactly one end and
cannot admit a proper \holo mapping onto a Riemann surface. In
fact, the following elementary observations suggest that \cpx
manifolds with the Bochner--Hartogs property are very different
topologically from those admitting proper \holo mappings onto
Riemann surfaces:

\begin{prop}\label{BH con boundary prop}
Let $X$ be a \con non\cpt \cpx manifold.
\begin{enumerate}

\item[(a)] Assume that $H^1_c(X,\ol)=0$. Then $e(X)=1$. In fact,
if $\Omega$ is any nonempty domain in~$X$ for which each \concomp
of the complement $X\sm\Omega$ is non\cptns, then $e(\Omega)=1$.
In particular, if $\Omega$ is a \rel \cpt domain in $X$ and
$X\sm\Omega$ is \conns, then $\partial\Omega$ is \conns. Moreover,
every \cpt orientable \cinf hypersurface in~$X$ is the boundary of
some smooth \rel \cpt domain in~$X$.

\item[(b)] If $X$ admits a surjective proper \cont open mapping
onto an orientable topological surface that is not simply \con
(for example, if $X$ admits a proper \holo mapping onto a Riemann
surface other than the disk or the plane), then there exists a
\cinf \rel \cpt domain~$\Omega$ in~$X$ \st $X\sm\Omega$ is \con
but $\partial\Omega$ is \emph{not} \con (and $e(\Omega)>1$). In
particular, $H^1_c(X,\ol)\neq 0$.

\end{enumerate}
\end{prop}
\begin{pf}
For the proof of (a), let us assume that $H^1_c(X,\ol)=0$. As
argued in the introduction, we must then have $e(X)=1$.  Next, we
show that any \cpt orientable \cinf hypersurface~$M$ in~$X$ is the
boundary of some \rel \cpt \cinf domain in~$X$. For we may fix a
\rel \cpt \con \nbdns~$U$ of $M$ in~$X$ \st $U\sm M$ has exactly
two \concompsns, $U_0$ and $U_1$. We may also fix a \rel \cpt
\nbdns~$V$ of~$M$ in~$U$ and a \cinf \fnns~$\lambda$ on~$X\sm M$
\st $\supp\lambda\Subset U$, $\lambda\equiv 0$ on $U_0\cap V$, and
$\lambda\equiv 1$ on $U_1\cap V$. Hence $\dbar\lambda$ extends to
a $\dbar$-closed \cinf $(0,1)$-form~$\alpha$ on~$X$ with \cpt
support in $U\sm M$, and since $H^1_c(X,\ol)=0$, we have
$\alpha=\dbar\beta$ for some \cinf compactly supported
\fnns~$\beta$ on~$X$. The difference $f\equiv\lambda-\beta$ is
then a \holo \fn on~$X\sm M$ that vanishes on some nonempty open
subset. If $X\sm M$ is \conns, then $f\equiv 0$ on the entire
set~$X\sm M$, and in particular, the restriction
$\beta\restriction_V$ is a \cinf \fn that is equal to~$1$ on
$U_1\cap V$, $0$ on $U_0\cap V$. Since $M=V\cap\partial
U_0=V\cap\partial U_1$, we have arrived at a contradiction. Thus
$X\sm M$ cannot be \conns, and hence $X\sm M$ must have exactly
two \concompsns, one containing $U_0$ and the other
containing~$U_1$. Since $e(X)=1$, one of these \concomps must be a
\rel \cpt \cinf domain with boundary~$M$ in~$X$. It follows that
in particular, the boundary of any \rel \cpt \cinf domain in~$X$
with \con complement must be \conns.

Next, suppose $\Omega$ is an arbitrary nonempty domain for which
$X\sm\Omega$ has no \cpt \compsns. If $e(\Omega)>1$, then part~(b)
of Lemma~\ref{basic ends lem} provides a \cinf \rel \cpt
domain~$\Theta$ in~$\Omega$ \st $\Omega\sm\Theta$ has no \cpt
\comps and $\partial\Theta$ is not \conns, and hence part~(e)
implies that $X\sm\Theta$ has no \cpt \compsns; i.e., $X\sm\Theta$
is \conns. However, as shown above, any smooth \rel \cpt domain
in~$X$ with \con complement must have \con boundary. Thus we have
arrived at a contradiction, and hence $\Omega$ must have only one
end. In particular, if $\Omega\Subset X$ (and $X\sm\Omega$ is
\conns), then by part~(c) of Lemma~\ref{basic ends lem},
$\partial\Omega$ must be \conns.

Part~(b) follows immediately from part~(f) of Lemma~\ref{basic
ends lem}.
\end{pf}

\section{Bounded geometry}\label{Bounded geometry sect}

In this section, we recall the definition of bounded geometry and
we fix some conventions. Let $X$ be a \cpx manifold with almost
\cpx structure $J\colon TX\to TX$. By a \emph{Hermitian metric}
on~$X$, we will mean a Riemannian metric~$g$ on $X$ \st
$g(Ju,Jv)=g(u,v)$ for every choice of real tangent vectors $u,v\in
T_pX$ with $p\in X$. We call $(X,g)$ a \emph{Hermitian manifold}.
We will also denote by $g$ the \cpx bilinear extension of $g$ to
the complexified tangent space $(TX)_{\C}$. The corresponding real
$(1,1)$-form~$\omega$ is given by $(u,v)\mapsto\omega(u,v)\equiv
g(Ju,v)$. The corresponding Hermitian metric (in the sense of a
smoothly varying family of Hermitian inner products) in the \holo
tangent bundle $T^{1,0}X$ is given by \((u,v)\mapsto g(u,\bar
v)\). Observe that with this convention, under the \holo vector
bundle isomorphism $(TX,J)\overset\cong\to T^{1,0}X$ given by
$u\mapsto \frac 12(u-iJu)$, the pullback of this Hermitian metric
to~$(TX,J)$ is given by \((u,v)\mapsto\tfrac 12g(u,v)-\tfrac
i2\omega(u,v)\). In a slight abuse of notation, we will also
denote the induced Hermitian metric in $T^{1,0}X$, as well as the
induced Hermitian metric in
$\Lambda^r(TX)_{\C}\otimes\Lambda^s(T^*X)_{\C}$, by~$g$. The
corresponding Laplacians are given by:
\begin{align*}
\lap&=\lap_d\equiv -(dd^*+d^*d),\\
\lap_{\dbar}&=-(\dbar\dbar^*+\dbar^*\dbar),\\
\lap_{\partial}&=-(\partial\partial^*+\partial^*\partial).
\end{align*}
If $(X,g,\omega)$ is \emph{K\"ahler}, i.e.,~$d\omega=0$, then
$\lap=2\lap_{\dbar}=2\lap_{\partial}$.

\begin{defn}\label{Bdd geom along set defn}
For $S\subset X$ and $k$ a nonnegative integer, we will say that a
Hermitian manifold~$(X,g)$ of dimension~$n$ {\it has bounded
geometry of order $k$ along $S$} if for some constant $C>0$ and
for every point $p\in S$, there is a biholomorphism $\Psi$ of the
unit ball $B\equiv B_{g_{\C^n}}(0;1)\subset\C^n$ onto a \nbd of
$p$ in $X$ \st $\Psi(0)=p$ and \st on $B$,
\[
C\inv g_{\C^n}\leq\Psi^*g\leq Cg_{\C^n}
\quad\text{and}\quad|D^m\Psi^*g|\leq C\text{ for }m=0,1,2,\dots,k.
\]
\end{defn}

\section{Green's \fns and harmonic projections}\label{greens fn sect}

In this section we recall some terminology and facts from
potential theory (a more detailed outline is provided in~\cite{NR-Structure theorems}). 
We will also see that the Bochner--Hartogs
property holds for a \con non\cpt complete K\"ahler manifold with
exactly one end and no nontrivial $L^2$ \holo $1$-forms.

A \con non\cpt oriented Riemannian manifold~$(M,g)$ is called
\textit{hyperbolic} if there exists a positive symmetric Green's
\fnns~$G(x,y)$ on~$M$; otherwise, $M$ is called \emph{parabolic}.
Equivalently, $M$ is hyperbolic if given a \rel \cpt \cinf
domain~$\Omega $ for which no \concomp of $M\sm\Omega$ is \cptns,
there is a \concompns~$E$ of $M\sm\overline\Omega$ and a (unique)
greatest \cinf \fn $u_E:\overline E\to[0,1)$ \st $u_E$ is \harm on
$E$, $u_E=0$ on $\partial E$, and $\sup_Eu_E=1$ (see, for example, Theorem~3 of
\cite{Glasner-Katz Function-theoretic degeneracy}). We will also call $E$, and any end containing~$E$, a \emph{hyperbolic} end. 
An end
that is not hyperbolic is called \emph{parabolic}, and 
we set $u_E\equiv 0$ for any parabolic end \compns~$E$ of $M\sm\Omega$.  We call the \fn $u:M\sm\Omega\to[0,1)$
defined by $u\restrict{\overline E}=u_E$ for each \concompns~$E$
of $M\sm\overline\Omega$, the \emph{\harm measure of the ideal
boundary of}~$M$ \emph{with respect to} $M\sm\overline\Omega$. A
sequence $\seq x\nu$ in $M$ with $x_\nu\to\infty$ and
$G(\cdot,x_\nu)\to 0$ (equivalently, $u(x_\nu)\to 1$) is called a
\emph{regular sequence}. Such a sequence always exists (for $M$
hyperbolic). A sequence $\seq x\nu$ tending to infinity with
$\liminf_{\nu\to\infty}G(\cdot,x_\nu)>0$ (i.e.,
$\limsup_{\nu\to\infty}u(x_\nu)<1$ or equivalently, $\seq x\nu$
has no regular subsequences) is called an \emph{irregular
sequence}. Clearly, every sequence tending to infinity that is not
regular admits an irregular subsequence. We say that an end $E$ of
$M$ is
\emph{regular} (\emph{irregular}) if every sequence in~$E$ tending to
infinity in~$M$ is regular (respectively, there exists an irregular
sequence in~$E$). Another characterization of hyperbolicity is that $M$
is hyperbolic if and only if $M$ admits a nonconstant negative
\cont sub\harm \fnns~$\vphi$. In fact, if $\seq x\nu$ is a
sequence in $M$ with $x_\nu\to\infty$ and $\vphi(x_\nu)\to 0$,
then $\seq x\nu$ is a regular sequence.

We recall that the \emph{energy} (or \emph{Dirichlet integral}) of
a suitable \fnns~$\vphi$ (for example, a \fn with first-order
distributional derivatives) on a Riemannian manifold~$M$ is given
by $\int_M|\nabla\vphi|^2\,dV$. To any \cinf compactly supported $\dbar$-closed $(0,1)$-form~$\alpha$ on a \con non\cpt hyperbolic
complete K\"ahler manifold $X$, we may associate a bounded finite-energy (i.e.,
Dirichlet-finite) \plh \fn on $X\sm\supp\alpha$ that vanishes
at infinity along any regular sequence:
\begin{lem}[see, for example, Lemma~1.1 of \cite{NR-BH Regular hyperbolic
Kahler}]\label{Cpt support dbar class gives plh fn lem} Let $X$ be
a \con non\cpt complete hyperbolic K\"ahler manifold, and let
$\alpha$ be a \cinf compactly supported form of type~$(0,1)$
on~$X$ with $\dbar\alpha=0$. Then there exist a closed and
coclosed $L^2$ harmonic form $\gamma$ of type $(0,1)$ and a \cinf
bounded \fn $\beta\colon X\to\C$ with finite energy \st
$\gamma=\alpha-\dbar\beta$ and $\beta(x_\nu)\to 0$ for every
regular sequence $\seq x\nu $ in $X$.
\end{lem}
\begin{rmks} 1. In particular, $\bar\gamma$
is a \holo $1$-form on $X$, and $\beta$ is \plh on the complement
of the support of~$\alpha$.

2. Under certain conditions, the leaves of the foliation
determined by $\bar\gamma$ outside a large \cpt subset of~$X$ are
\cptns, and one gets a proper \holo mapping onto a Riemann surface.

3. According to Lemma~\ref{holo vanish along regular sequence not
hyp lem} below (which is a modification of an observation due to
J.~Wang), if $\beta$ is \holo on some hyperbolic end, then $\beta$
vanishes on that end.
\end{rmks}

\begin{lem}[cf.~Lemma~1.3 of \cite{NR-BH Regular hyperbolic
Kahler}]\label{holo vanish along regular sequence not hyp lem} Let
$X$ be a \con non\cpt complete (hyperbolic) K\"ahler manifold, and
let $E$ be a hyperbolic end of~$X$. If $f$ is a bounded
\holo \fn on~$E$ and $f(x_\nu)\to 0$ for every regular
sequence~$\seq x\nu$ for $X$ in~$E$, then $f\equiv 0$ on~$E$.
\end{lem}
\begin{pf}
We may fix a nonempty smooth domain~$\Omega$ \st
$\partial E\subset\Omega\Subset X$ and $X\sm\Omega$ has no \cpt
\concompsns. In particular, some \compns~$E_0$ of $E\sm\overline\Omega$ is a hyperbolic end of~$X$. The \hmib $X$ \wrt $X\sm\overline\Omega$ is a
nonconstant \fn $u\colon X\sm\Omega\to[0,1)$. By replacing~$f$
with the product of~$f$ and a sufficiently small nonzero constant,
we may assume that $|f|<1$ and hence for each $\epsilon>0$,
$u+\epsilon\log|f|<0$ on $E\cap\partial\Omega$. Thus we get a
nonnegative bounded \cont sub\harm \fn $\vphi_\epsilon$ on $X$ by setting
$\vphi_\epsilon\equiv 0$ on $X\sm E_0$ and
$\vphi_\epsilon\equiv\max(0,u+\epsilon\log|f|)$ on~$E_0$.
If $f(p)\neq 0$ at some point $p\in E_0$, then
$\vphi_\epsilon(p)>0$ for $\epsilon$ sufficiently small. However,
any sequence $\seq x\nu$ in~$E_0$ with $\vphi_\epsilon(x_\nu)\to
m\equiv\sup\vphi_\epsilon>0$ must be a regular sequence and must
therefore satisfy \(u(x_\nu)+\epsilon\log|f(x_\nu)|\to-\infty\),
which contradicts the choice of $\seq x\nu$. Thus $f$ vanishes on
$E_0$ and therefore, on~$E$.
\end{pf}

The above considerations lead to the following observation
(cf.~Proposition~4.4 of~\cite{NR Weak Lefschetz}):
\begin{thm}\label{Bochner Hartogs no L2 forms thm}
Let $X$ be a \con non\cpt hyperbolic complete K\"ahler manifold
with no nontrivial $L^2$ \holo $1$-forms.
\begin{enumerate}

\item[(a)] For every compactly supported $\dbar$-closed \cinf
form~$\alpha$ of type~$(0,1)$ on~$X$, there exists a bounded \cinf
\fnns~$\beta$ with finite energy on~$X$ \st $\dbar\beta=\alpha$
on~$X$ and $\beta$ vanishes on every hyperbolic end $E$ of~$X$
that is contained in $X\sm\supp\alpha$.

\item[(b)] In any ends decomposition $X\sm K=E_1\cup\cdots\cup
E_m$, exactly one of the ends, say~$E_1$, is hyperbolic, and
moreover, every \holo \fn on~$E_1$ admits a (unique) extension to
a \holo \fn on~$X$.

\item[(c)] If $e(X)=1$ (equivalently, every end of~$X$ is
hyperbolic), then $H^1_c(X,\ol)=0$.

\end{enumerate}
\end{thm}
\begin{pf}
Given a compactly supported $\dbar$-closed \cinf form~$\alpha$ of
type~$(0,1)$ on~$X$, Lemma~\ref{Cpt support dbar class gives plh
fn lem} provides a bounded \cinf \fn $\beta$ with finite energy
\st $\dbar\beta=\alpha$ and $\beta(x_\nu)\to 0$ for every regular
sequence $\seq x\nu$ in $X$ (by hypothesis, the $L^2$ \holo
$1$-form~$\bar\gamma$ provided by the lemma must be trivial). In
particular, $\beta$ is \holo on $X\sm\supp\alpha$, and
Lemma~\ref{holo vanish along regular sequence not hyp lem} implies
that $\beta$ must vanish on every hyperbolic end of~$X$ contained
in~$X\sm\supp\alpha$. Thus part~(a) is proved.

For the proof of part~(b), suppose $X\sm K=E_1\cup\cdots\cup E_m$
is an ends decomposition. Then at least one of the ends,
say~$E_1$, must be hyperbolic. Given a \fn $f\in\ol(X\sm K)$, we
may fix a \rel \cpt \nbdns~$U$ of $K$ in~$X$ and a \cinf
\fnns~$\lambda$ on~$X$ \st $\lambda\equiv f$ on~$X\sm U$. Applying
part~(a) to the $(0,1)$-form $\alpha\equiv\dbar\lambda$, we get a
\cinf \fnns~$\beta$ \st $\dbar\beta=\alpha$ on~$X$ and
$\beta\equiv 0$ on any hyperbolic end contained in~$X\sm U$. If
$E_j$ is a hyperbolic end (for example, if $j=1$), then $E_j\sm U$
must contain a hyperbolic end~$E$ of~$X$, and the \holo \fn
$h\equiv\lambda-\beta$ on~$X$ must agree with~$f$ on~$E$ and
therefore, on~$E_j$. Thus we get a \holo \fn on~$X$ that agrees
with~$f$ on every $E_j$ which is hyperbolic. Taking $f$ to be a
locally constant \fn on $X\sm K$ with distinct values on the
\comps $E_1,\dots,E_m$, we see that in fact, $E_j$ must be a
parabolic end for $j=2,\dots,m$.

Part~(c) follows immediately from parts~(a) and (b).
\end{pf}

We close this section with a preliminary step toward the proof of
Theorem~\ref{BHD bounded geom hyperbolic thm}:
\begin{lem}\label{BHD hyperbolic with plh bounded gradient etc lem}
Suppose $(X,g)$ is a \con non\cpt hyperbolic complete K\"ahler
manifold with bounded geometry of order~$0$, $e(X)=1$, and there
exists a real-valued \plh \fnns~$\rho$ with bounded gradient and
infinite energy on~$X$. Then $X$ admits a proper \holo mapping
onto a Riemann surface if and only if $H^1_c(X,\ol)\neq 0$.
\end{lem}
\begin{pf}
Given a compactly supported $\dbar$-closed \cinf form~$\alpha$ of
type~$(0,1)$ on~$X$, Lemma~\ref{Cpt support dbar class gives plh
fn lem} provides a closed and coclosed $L^2$ harmonic form
$\gamma$ of type $(0,1)$ and a \cinf bounded \fn $\beta\colon
X\to\C$ with finite energy \st $\gamma=\alpha-\dbar\beta$ and
$\beta(x_\nu)\to 0$ for every regular sequence $\seq x\nu $ in
$X$. If $\gamma\equiv 0$, then $\dbar\beta=\alpha$ and
Lemma~\ref{holo vanish along regular sequence not hyp lem} implies
that~$\beta$ vanishes on the complement of some \cpt set. If
$\gamma$ is nontrivial, then the $L^2$ \holo $1$-form
$\theta_1\equiv\bar\gamma$ and the bounded \holo $1$-form
$\theta_2\equiv\partial\rho$, which is not in~$L^2$, must be
linearly independent.  Theorem~0.1 and Theorem~0.2 of~\cite{NR L2
Castelnuovo} then provide a proper \holo mapping of~$X$ onto a
Riemann surface.
\end{pf}
\begin{rmk}
The proofs of Lemma~1.1 of~\cite{NR L2 Castelnuovo} and
Theorem~0.1 of~\cite{NR L2 Castelnuovo} (the latter fact was
applied above and relies on the former) contain a minor mistake in
their application of continuity of intersections (see \cite{Stein} 
or \cite{Tworzewski-Winiarski Cont of intersect} or Theorem~4.23
in \cite{ABCKT}). In each of these proofs, one has a sequence of
levels $\seq L\nu$ of a \holo mapping $f\colon X\to\mathbb P^1$
and a sequence of points $\seq x\nu$ \st $x_\nu\in L_\nu$ for
each~$\nu$ and $x_\nu\to p$. For $L$ the level of~$f$ through~$p$,
by continuity of intersections, $\seq L\nu$ converges to~$L$
relative to the ambient manifold $X\sm[f\inv(f(p))\sm L]$, but
contrary to what was stated in these proofs, a~priori, this
convergence need not hold relative to~$X$. Aside from this small
misstatement, the proofs are correct and no further changes are
needed.
\end{rmk}

\section{Quasi-Dirichlet-finite pluriharmonic
functions}\label{quasidirichlet finite fns sect}

The following is the main advantage of working with
\emph{irregular} hyperbolic manifolds:
\begin{lem}[Nakai]\label{quasidirichlet finite irreg hyp exists lem}
Let $(M,g)$ be a \con non\cpt irregular hyperbolic oriented
complete Riemannian manifold, let $\seq qk$ be an irregular
sequence, let $G(\cdot,\cdot)$ be the Green's \fnns, and let
\(\rho_k\equiv G(\cdot,q_k)\colon M\to(0,\infty]\) for each~$k$.
Then some subsequence of $\seq\rho k$ converges uniformly on \cpt
subsets of~$M$ to a \fnns~$\rho$. Moreover, any such limit
\fnns~$\rho$ has the following properties:
\begin{enumerate}

\item[(i)] The \fnns~$\rho$ is positive and \harmns;

\item[(ii)] $\int_M|\nabla\rho|^2\,dV_g=\infty$;

\item[(iii)] \(\int_{\rho\inv([a,b])}|\nabla\rho|_g^2\,dV_g\leq
b-a\) for all $a$ and $b$ with $0\leq a<b$ (in particular, $\rho$
is unbounded); and

\item[(iv)] If $\Omega$ is any smooth domain with \cpt boundary (i.e., either $\Omega$ is an end or $\Omega\Subset M$)
and at most finitely many terms of the sequence $\seq qk$ lie
in~$\Omega$, then
\[
\sup_\Omega\rho=\max_{\partial\Omega}\rho<\infty\qquad\text{and}\qquad
\int_\Omega|\nabla\rho|^2\,dV\leq\int_{\partial\Omega}\rho
\pdof\rho\nu\,d\sigma<\infty.
\]

\end{enumerate}
\end{lem}
\begin{rmk}
Following Nakai \cite{Nakai Green potential} and Sario and Nakai
\cite{SarioNakai}, a positive \fn $\vphi$ on a Riemannian manifold
$(M,g)$ is called \emph{quasi-Dirichlet-finite} if there is a
positive constant~$C$ \st
\[
\int_{\vphi\inv([0,b])}|\nabla\vphi|_g^2\,dV_g\leq Cb
\]
for every $b>0$. Nakai proved the existence of an Evans-type
quasi-Dirichlet-finite positive \harm \fn on an irregular Riemann
surface. His arguments, which involve the behavior of the Green's
\fn at the Royden boundary, carry over to a Riemannian manifold
and actually show that the constructed \fn has the slightly
stronger property appearing in the above lemma. One can instead
prove the lemma via Nakai's arguments simply by taking
$\rho=G(\cdot, q)$, where $G$ is the extension of the Green's \fn
to the Royden compactification and $q$ is a point in the Royden
boundary for which $\rho>0$ on $M$. The direct proof appearing
below is essentially this latter argument.
\end{rmk}
\begin{pf*}{Proof of Lemma~\ref{quasidirichlet finite irreg hyp exists
lem}} Fixing a sequence of nonempty smooth domains $\seq\Omega
m_{m=0}^\infty$ \st $M\sm\Omega_0$ has no \cpt \concompsns,
$\bigcup_{m=0}^\infty\Omega_m=M$, and
$\Omega_{m-1}\Subset\Omega_m$ for $m=1,2,3,\dots$, and letting
$G_m$ be the Green's \fn on $\Omega_m$ for each~$m$, we get
$G_m\nearrow G$. Given $m_0\in\Z_{>0}$, for each integer $m>m_0$
and each point $p\in\Omega_{m_0}$, the \cont
\fnns~$G_m(p,\cdot)\restrict{\overline\Omega_m\sm\Omega_{m_0}}$
vanishes on~$\partial\Omega_m$, and the \fn is positive on
$\partial\Omega_{m_0}$ and \harm on
$\Omega_m\sm\overline\Omega_{m_0}$. Thus
\[
G_m(p,\cdot)\restrict{\overline\Omega_m\sm\Omega_{m_0}}
\leq\max_{\partial\Omega_{m_0}}G_m(p,\cdot)\leq\max_{\partial\Omega_{m_0}}
G(p,\cdot).
\]
Passing to the limit we get
\[
G(p,\cdot)\leq\max_{\partial\Omega_{m_0}}G(p,\cdot)
\]
on $M\sm\Omega_{m_0}$ for each point $p\in\Omega_{m_0}$. Hence
\[
G\leq
A_{m_0}\equiv\max_{\overline\Omega_{m_0-1}\times\partial\Omega_{m_0}}G
\]
on $\overline\Omega_{m_0-1}\times(M\sm\Omega_{m_0})$. In
particular, $\rho_k=G(\cdot,q_k)\leq A_{m_0}$ on
$\overline\Omega_{m_0-1}$ for $k\gg 0$. Therefore, by replacing
$\seq qk$ with a suitable subsequence, we may assume that $\rho_k$
converges uniformly on \cpt subsets of~$M$ to a positive \harm
\fnns~$\rho$.

Suppose $0<a<b$. Given $k\in\Z_{>0}$, for $m\gg 0$ we have
$q_k\in\Omega_m$, and the \fn $\rho_k\ssp m\equiv
G_m(\cdot,q_k)\colon\overline\Omega_m\to[0,\infty]$ satisfies
\((\rho_k\ssp m)\inv((a,\infty])\Subset\Omega_m\). Hence if $r$ and $s$
are regular values of $\rho_k\ssp m\restrict{\Omega_m\sm\set{q_k}}$
with $a<r<s<b$, then
\begin{align*}
\int_{(\rho_k\ssp m)\inv((r,s))}|\nabla\rho_k\ssp m|^2\,dV&
=\int_{(\rho_k\ssp m)\inv(s)}\rho_k\ssp m\pdof{\rho_k\ssp m}{\nu}\,d\sigma
-\int_{(\rho_k\ssp m)\inv(r)}\rho_k\ssp m\pdof{\rho_k\ssp m}{\nu}\,d\sigma\\
&=\int_{(\rho_k\ssp m)\inv(s)}s\cdot\pdof{\rho_k\ssp m}\nu\,d\sigma
-\int_{(\rho_k\ssp m)\inv(r)}r\cdot\pdof{\rho_k\ssp m}\nu\,d\sigma\\
&=(r-s)\int_{\partial\Omega_m}\pdof{\rho_k\ssp m}\nu\,d\sigma\\
&=(s-r)\int_{\partial\Omega_m}(-1)\pdwrt\nu\left[G_m(\cdot,q_k)\right]\,d\sigma\\
&=s-r,
\end{align*}
where $\partial/\partial\nu$ is the normal derivative oriented
outward for the open sets $\Omega_m$, $(\rho_k\ssp m)\inv((0,s))$, and
$(\rho_k\ssp m)\inv((0,r))$.  Here we have normalized $G_m$ (and
similarly, all Green's \fnsns) so that $-\lap_{\text{distr.}}
G_m(\cdot,q)$ is the Dirac \fn at~$q$ for each
point~$q\in\Omega_m$. 
Letting $r\to a^+$ and $s\to b^-$, we get
\[
\int_{(\rho_k\ssp m)\inv((a,b))}|\nabla\rho_k\ssp m|^2\,dV=(b-a).
\]
Letting $\chi_A$ denote the characteristic \fn of each set
$A\subset M$, we have
\[
\lim_{m\to\infty}|\nabla\rho_k\ssp m|=|\nabla\rho_k|\text{ on
}M\sm\set{q_k}\qquad\text{and}\qquad
\liminf_{m\to\infty}\chi_{(\rho_k\ssp m)\inv((a,b))}\geq\chi_{\rho_k\inv((a,b))}.
\]
Hence Fatou's lemma gives
\[
\int_{\rho_k\inv((a,b))}|\nabla\rho_k|^2\,dV\leq(b-a).
\]
Similarly, letting $k\to\infty$, we get
\[
\int_{\rho\inv((a,b))}|\nabla\rho|^2\,dV\leq(b-a).
\]
Applying this inequality to $a'$ and $b'$ with $0<a'<b'$ and
letting $a'\to a^-$ and $b'\to b^+$, we get
\[
\int_{\rho\inv([a,b])}|\nabla\rho|^2\,dV\leq(b-a).
\]
Letting $a\to 0^+$ (and noting that $\rho>0$), we also get the
above inequality for $a=0$.

Assuming now that $\rho$ has finite energy, we will reason to a
contradiction. We may fix a
constant~$b>\sup_{\Omega_0}\rho$ that is a regular value of~$\rho$, of
$\rho_k\restrict{M\sm\set{q_k}}$ for all~$k$, and of $\rho_k\ssp m\restrict{\Omega_m\sm\set{q_k}}$ for all
$k$ and~$m$. Note that we have not yet shown that~$\rho$ is unbounded,
so we have not yet ruled out the possibility that $\rho\inv((0,b))=M$, and in particular, that $\rho\inv(b)=\emptyset$. Since $\rho_k\to\rho$ uniformly on \cpt subsets of~$M$ as
$k\to\infty$, and for each~$k$, $\rho_k\ssp m\to\rho_k$ uniformly on
\cpt subsets of $M\sm\set{q_k}$ as $m\to\infty$, we may fix a positive integer~$k_0$ and
a \str increasing sequence of positive integers $\seq mk$ \st
$q_k\in\Omega_{m_k}$ for each~$k$, 
$\rho_k\ssp{m_k}\leq\rho_k<b$ on
$\overline\Omega_0$ for each~$k\geq k_0$, and $\rho_k\ssp{m_k}\to\rho$
uniformly on \cpt sets as $k\to\infty$.
Letting $\vphi\equiv\min(\rho,b)$ and letting $\vphi_k\colon
M\to[0,b]$ be the Lipschitz \fn given by
\[
\vphi_k\equiv
\begin{cases}
{\min(\rho_k\ssp{m_k},b)}&{\text{ on }\overline\Omega_{m_k}}\\
{0}&{\text{ elsewhere}}
\end{cases}
\]
for each~$k$, we see that $\vphi_k\to\vphi$ uniformly on \cpt
subsets of~$M$ and $\nabla\vphi_k\to\nabla\vphi$ uniformly on \cpt
subsets of $M\sm\rho\inv(b)$. Moreover, for each $k$,
\[
\int_M|\nabla\vphi_k|^2\,dV
=\int_{(\rho_k\ssp{m_k})\inv((0,b))}|\nabla\rho_k\ssp{m_k}|^2\,dV=b.
\]
Applying weak compactness, we may assume that
$\set{\nabla\vphi_k}$ converges weakly in $L^2$ to a vector
field~$v$. But for each \cpt set $K\subset M\sm\rho\inv(b)$,
$(\nabla\vphi_k)\restrict K\to(\nabla\vphi)\restrict K$ uniformly,
and therefore in~$L^2$. Since $\rho\inv(b)$ is a set of
measure~$0$, we must have $v=\nabla\vphi$ (in~$L^2$).
Hence
\begin{align*}
\int_{\rho\inv((0,b))}|\nabla\rho|^2\,dV
=\langle\nabla\vphi,\nabla\rho\rangle\leftarrow\langle\nabla\vphi_k,\nabla\rho\rangle
&=\int_{\partial\Omega_{m_k}}\rho_k\ssp{m_k}\pdof
\rho\nu\,d\sigma\\
&\qquad\qquad\qquad+\int_{(\rho_k\ssp{m_k})\inv(b)}\rho_k\ssp{m_k}\pdof
\rho\nu\,d\sigma\\
&=0-b\int_{\partial\left((\rho_k\ssp{m_k})\inv((b,\infty])\right)}
\pdof\rho\nu\,d\sigma=0.
\end{align*}
It follows that $\rho\equiv a$ for some constant~$a$ (in
particular, $0<a<b$). Letting $u$ be the \hmib $M$ with respect to
$M\sm\overline\Omega_0$ and letting $\psi\colon M\to[0,1)$ be the
Dirichlet-finite locally Lipschitz \fn on $M$ obtained by extending $u$
by~$0$, we get
\begin{align*}
0=\langle 0,\nabla\psi\rangle
\leftarrow\langle\nabla\vphi_k,\nabla\psi\rangle
&=\int_{\partial\Omega_{m_k}}\rho_k\ssp{m_k}\pdof
u\nu\,d\sigma-\int_{\partial\Omega_0}\rho_k\ssp{m_k}\pdof
u\nu\,d\sigma\\
&\qquad\qquad\qquad\qquad\qquad+\int_{(\rho_k\ssp{m_k})\inv(b)}\rho_k\ssp{m_k}\pdof u\nu\,d\sigma\\
&=0-\int_{\partial\Omega_0}\rho_k\ssp{m_k}\pdof
u\nu\,d\sigma-b\int_{\partial\left((\rho_k\ssp{m_k})\inv((b,\infty])\right)}
\pdof u\nu\,d\sigma\\
&=-\int_{\partial\Omega_0}\rho_k\ssp{m_k}\pdof
u\nu\,d\sigma\to-a\int_{\partial\Omega_0}\pdof u\nu\,d\sigma<0.
\end{align*}
Thus we have arrived at a contradiction, and hence $\rho$ must
have infinite energy.

Finally, given a smooth domain~$\Omega$ as in~(iv), for each $k\gg
0$, we have $q_k\notin\overline\Omega$. For $m\gg 0$, we have
$q_k\in\Omega_m$ and $\partial\Omega\subset\Omega_m$. Since
$\rho_k\ssp m$ is then \cont on $\overline{\Omega\cap\Omega_m}$,
\harm on $\Omega\cap\Omega_m$, and zero on $\partial\Omega_m$, we
also have
\[
\sup_{\Omega\cap\Omega_m}\rho_k\ssp m=\max_{\partial\Omega}\rho_k\ssp m\qquad\text{and}\qquad
\int_{\Omega\cap\Omega_m}|\nabla\rho_k\ssp m|^2\,dV
=\int_{\partial\Omega}\rho_k\ssp m\pdof{\rho_k\ssp m}\nu\,d\sigma.
\]
Letting $m\to\infty$, and then letting $k\to\infty$, we get the
required properties of~$\rho$ on~$\Omega$.
\end{pf*}

We will also use the following analogue of a theorem of Sullivan
(see~\cite{Sullivan} and Theorem~2.1 of~\cite{NR-Structure
theorems}):
\begin{lem}\label{quasidirichlet harm bdd grad lem}
Let $(M,g)$ be a \con non\cpt oriented complete Riemannian
manifold, $E$ an end of~$M$, and $h$ a positive \cinf \fn on~$M$.
Assume that:
\begin{enumerate}

\item[(i)]  There exist positive constants~$K$, $R_0$, and
$\delta$ \st \(\text{\rm Ric}_g\geq-Kg\) on~$E$ and
\(\Vol\,(B(x;R_0))\geq\delta\) for every point $x\in E$;

\item[(ii)] The restriction $h\restrict E$ is \harmns; and

\item[(iii)] For some positive constant $C$, $\int_{E\cap
h\inv([a,b])}|\nabla h|^2\,dV\leq C(b-a)+C$ for all $a$ and $b$
with $0\leq a<b$.

\end{enumerate}
Then $|\nabla h|$ is bounded on~$E$, and for each point $p\in M$,
\[
\int_{E\cap B(p;R)}|\nabla h|^2\,dV=O(R)\text{ as }R\to\infty.
\]
\end{lem}
\begin{pf*}{Sketch of the proof}
We may fix a nonempty \cpt set $A\supset\partial E$. As in the
proof of Theorem~2.1 of~\cite{NR-Structure theorems}, setting
$\vphi\equiv|\nabla h|^2$, we get a positive constant $C_1$ \st
for each point $x_0\in E$ with $\text{dist}(x_0,A)>R_1\equiv
4R_0$,
\[
\sup_{B(x_0;R_0)}\vphi\leq C_1\int_{B(x_0;2R_0)}\vphi\,dV.
\]
For $a\equiv\inf_{B(x_0;2R_0)}h$ and $b\equiv\sup_{B(x_0;2R_0)}h$,
we have on the one hand,
\[
\int_{B(x_0;2R_0)}\vphi\,dV\leq\int_{E\cap
h\inv([a,b])}\vphi\,dV\leq C(b-a)+C.
\]
On the other hand, $b-a\leq\sup_{B(x_0,R_1)}|\nabla h|R_1$.
Combining the above, we see that if $C_2>1$ is a sufficiently
large positive constant that is, in particular, greater than the
supremum of $|\nabla h|$ on the $2R_1$-\nbd of~$A$, then for each
point $x_0\in E$ at which $|\nabla h(x_0)|>C_2$, we have
\[
\sup_{B(x_0;R_0)}|\nabla h|^2<C_2\sup_{B(x_0;R_1)}|\nabla h|.
\]
Fixing constants $C_3>C_2$ and $\epsilon>0$ so that
$C_3^{1-\epsilon}>C_2$, we see that if $|\nabla h(x_0)|>C_3$, then
there exists a point $x_1\in B(x_0;R_1)$ \st
\[
(1+\epsilon)\log|\nabla h(x_0)|\leq\log|\nabla h(x_1)|.
\]

Assuming now that $|\nabla h|$ is unbounded on~$E$, we will reason
to a contradiction.  Fixing a point $x_0\in E$ at which $|\nabla
h(x_0)|>C_3$ and applying the above inequality inductively, we get
a sequence $\seq xm$ in $E$ \st $\text{dist}(x_m,x_{m-1})<R_1$ and
\[
(1+\epsilon)\log|\nabla h(x_{m-1})|\leq\log|\nabla h(x_m)|
\]
for $m=1,2,3,\dots$; that is, $\set{|\nabla h(x_m)|}$ has
super-exponential growth. However, the local version of Yau's
Harnack inequality (see \cite{Cheng-Yau Diff eqs mflds}) provides
a constant $C_4>0$ \st 
\[
|\nabla h(x)|\leq C_4h(x)\qquad\text{and}\qquad h(x)\leq C_4 h(p)
\]
for all points $x,p\in M$ with $\text{dist}(p,A)>2R_1$ and
$\text{dist}(x,p)<R_1$, so $\set{|\nabla h(x_m)|}$ has at most
exponential growth. Thus we have arrived at a contradiction, and
hence $|\nabla h|$ must be bounded on~$E$.

Finally, by redefining $h$ outside a \nbd of~$\overline E$, we may
assume without loss of generality that $|\nabla h|$ is bounded
on~$M$. Fixing a point $p\in M$, we see that for $R>0$,
$a\equiv\inf_{B(p;R)}h$, and $b\equiv\sup_{B(p;R)}h$, we have
\[
\int_{E\cap B(p;R)}|\nabla h|^2\,dV\leq\int_{E\cap
h\inv([a,b])}|\nabla h|^2\,dV\leq C(b-a)+C\leq C\cdot\sup|\nabla
h|\cdot 2R+C.
\]
Therefore,
\[
\int_{E\cap B(p;R)}|\nabla h|^2\,dV=O(R)\text{ as }R\to\infty.
\]
\end{pf*}

Applying the above in the K\"ahler setting, we get the following:
\begin{prop}\label{plh fn infinite energy exist kahler prop}
Let $(X,g)$ be a \con non\cpt complete K\"ahler manifold, let $E$ be
an irregular hyperbolic end along which $X$ has bounded geometry
of order~$2$ (or for which there exist positive constants~$K$, $R_0$, and
$\delta$ \st \(\text{\rm Ric}_g\geq-Kg\) on~$E$ and
\(\Vol\,(B(x;R_0))\geq\delta\) for every point $x\in E$), let $\seq qk$ be an irregular sequence in~$E$, let
$G(\cdot,\cdot)$ be the Green's \fn on~$X$, and let \(\rho_k\equiv
G(\cdot,q_k)\colon X\to(0,\infty]\) for each~$k$. Then some
subsequence of $\seq\rho k$ converges uniformly on \cpt subsets
of~$X$ to a \fnns~$\rho$. Moreover, any such limit \fnns~$\rho$
has the following properties:
\begin{enumerate}

\item[(i)] The \fnns~$\rho$ is positive and \plhns;

\item[(ii)] $\int_E|\nabla\rho|^2\,dV=\infty>\int_{X\sm
E}|\nabla\rho|^2\,dV$;

\item[(iii)] \(\int_{\rho\inv([a,b])}|\nabla\rho|^2\,dV\leq b-a\)
for all $a$ and $b$ with $0\leq a<b$ (in particular, $\rho$ is
unbounded on~$E$); 

\item[(iv)] If $\Omega$ is any smooth domain with \cpt boundary (i.e., either $\Omega$ is an end or $\Omega\Subset X$)
and at most finitely many terms of the sequence $\seq qk$ lie
in~$\Omega$, then
\[
\sup_\Omega\rho=\max_{\partial\Omega}\rho<\infty\qquad\text{and}\qquad
\int_\Omega|\nabla\rho|^2\,dV\leq\int_{\partial\Omega}\rho
\pdof\rho\nu\,d\sigma<\infty;
\]
and

\item[(v)] $|\nabla\rho|$ is bounded.

\end{enumerate}
\end{prop}
\begin{pf}
By Lemma~\ref{quasidirichlet finite irreg hyp exists lem}, some subsequence of $\seq\rho k$ converges uniformly on \cpt sets, and the limit~$\rho$ of any such subsequence is positive and \harm and satisfies (ii)--(iv).
Lemma~\ref{quasidirichlet harm bdd grad lem} implies that
$|\nabla\rho|$ is bounded on~$E$ and for $p\in X$,
\(\int_{B(p;R)}|\nabla\rho|^2\,dV=O(R)\) (hence $\int_{B(p;R)}|\nabla\rho|^2\,dV=o(R^2)$) as
\(R\to\infty\). By an observation of Gromov~\cite{Gro-Kahler
hyperbolicity} and of Li~\cite{Li Structure complete Kahler} (see
Corollary~2.5 of~\cite{NR-Structure theorems}), $\rho$ is \plhns.
\end{pf}

\section{Proof of the main result and some related results}\label{proof main thm sect}

This section contains the proof of Theorem~\ref{BHD bounded geom
hyperbolic thm}.  We also consider some related results.
\begin{pf*}{Proof of Theorem~\ref{BHD bounded geom hyperbolic
thm}} Let $X$ be a \con non\cpt hyperbolic complete K\"ahler
manifold with \bdd geometry of order two, and assume that $X$ has
exactly one end. By the main result of~\cite{NR-BH Regular
hyperbolic Kahler}, the Bochner--Hartogs dichotomy holds for $X$
regular hyperbolic. If $X$ is irregular hyperbolic, then
Proposition~\ref{plh fn infinite energy exist kahler prop}
provides a (quasi-Dirichlet-finite) positive \plh \fn~$\rho$
on~$X$ with infinite energy and bounded gradient, and hence
Lemma~\ref{BHD hyperbolic with plh bounded gradient etc lem} gives
the claim.
\end{pf*}

The above arguments together with those appearing in
\cite{NR-Structure theorems}, \cite{NR Filtered ends}, and
\cite{NR L2 Castelnuovo} give results for multi-ended complete
K\"ahler manifolds. To see this, we first recall some terminology and facts.

\begin{defn}\label{ends filtered ends def}
Let $M$ be a \con manifold. Following Geoghegan \cite{Geoghegan} (see also
Kropholler and Roller \cite{KroR}), for $\Upsilon\colon\widetilde
M\to M$ the universal covering of $M$, elements of the set
\[
\lim_{\leftarrow}\pi_0 [\Upsilon\inv(M\sm K)],
\]
where the limit is taken as $K$ ranges over the \cpt subsets of
$M$ (or the \cpt subsets of $M$ for which the complement $M\setminus K$
has no \rel \cpt \compsns) will be called {\it filtered ends}. The
number of filtered ends of $M$ will be denoted by~$\tilde e(M)$.
\end{defn}

\begin{lem}\label{Behavior of filtered ends for covering lem}
Let $M$ be a \con non\cpt topological manifold.
\begin{enumerate}

\item[(a)] We have $\tilde e(M)\geq e(M)$. In fact for any $k\in
\N$, we have $\tilde e(M)\geq k$ if and only if there exists an
ends decomposition $M\setminus K=E_1\cup \cdots \cup E_m$ \st
\[
\sum _{j=1}^m[\pi _1(M):\Gamma _j]\geq k,
\]
where $\Gamma_j\equiv\image{\pi_1(E_j)\to\pi_1(M)}$ for
$j=1,\dots,m$.

\item[(b)] If $\Upsilon\colon\widehat M\to M$ is a \con covering
space, $E$ is an end of~$\what M$, and
$E_0\equiv\Upsilon(E)\varsubsetneqq M$, then
\begin{enumerate}

\item[(i)] $E_0$ is an end of~$M$;

\item[(ii)] $\partial E_0=\Upsilon(\partial E)\sm E_0$;

\item[(iii)] $\overline E\cap\Upsilon\inv(\partial
E_0)=(\partial E)\sm\Upsilon\inv(E_0)$;

\item[(iv)] The mapping $\Upsilon\restrict{\overline
E}\colon\overline E\to\overline E_0$ is proper and surjective; and

\item[(v)] If $F_0\subset E_0\sm\Upsilon(\partial E)$ is an end
of~$M$ and $F\equiv E\cap\Upsilon\inv(F_0)$, then
$\Upsilon\restrict F\colon F\to F_0$ is a finite covering and each
\concomp of $F$ is an end of~$\what M$.

\end{enumerate}

\item[(c)]  If $\Upsilon\colon\widehat M\to M$ is a \con covering
space, then $\tilde e(\widehat M)\leq\tilde e(M)$, with equality
holding if the covering is finite.

\end{enumerate}
\end{lem}
\begin{pf}
For any nonempty domain~$U$ in~$M$, the index of $\image{\pione
U\to\pione M}$ is equal to the number of \concomps of the lifting
of $U$ to the universal covering of~$M$, so part~(a) holds.

For the proof of part~(b), observe that $E_0$ is a domain in~$M$,
$\partial E_0\neq\emptyset$, $\Upsilon(\overline E)\subset\overline
E_0$, and therefore, $\overline E\cap\Upsilon\inv(\partial
E_0)=(\partial E)\sm\Upsilon\inv(E_0)$. Given a point $p\in M$, we may fix domains $U$
and $V$ in $M$ \st $p\in U\Subset V$, $U\cap\partial E_0\neq\emptyset$,
and the image of $\pione V$ in $\pione X$ is trivial (their existence 
is trivial if $p\in\partial E_0$, while for $p\notin\partial E_0$, we may take $U$ 
and $V$ to be sufficiently small \con \nbds of the image of an
injective path from~$p$ to a point in $\partial E_0$).
The \concomps of $\what U\equiv\Upsilon\inv(U)$ then form a locally finite
collection of \rel \cpt domains in $\what M$, and those \comps
that meets~$\overline E$ must also meet the \cpt set $\partial E$, so only
finitely many \compsns, say $U_1,\dots,U_m$,
meet~$\overline E$. Thus for $U_0\equiv\bigcup_{i=1}^mU_i$,
we have $\what U\cap\overline E=U_0\cap\overline E\Subset\what M$, 
and it follows that the restriction
$\overline E\to\overline E_0$ is a proper mapping. In particular, this is a closed mapping, and hence $\Upsilon(\overline E)=\overline E_0$.
Furthermore, the boundary
\[
\partial E_0=\Upsilon(\overline E)\sm\Upsilon(E)=\Upsilon(\partial E)\sm E_0
\]
is \cpt and $\overline E_0$ is non\cpt (by properness), so $E_0$ must be an end of~$M$.

Finally, if $F_0\subset E_0\sm\Upsilon(\partial E)$ is an end
of~$M$, then each \concomp of $\Upsilon\inv(F_0)$ that meets $E$
must lie in~$E$. Thus the
restriction $F\equiv E\cap\Upsilon\inv(F_0)\to F_0$ is a covering
space.  Properness then implies that this restriction is actually
a finite covering, $\partial F\subset\overline
E\cap\Upsilon\inv(\partial F_0)$ is \cptns, and in particular,
each \concomp of~$F$ is an end of~$\what M$.

For the proof of part~(c), let $\what\Upsilon\colon\wtil M\to\what
M$ be the universal covering, and let $k\in\N$ with $\tilde
e(\what M)\geq k$. Then there exists an ends decomposition $\what
M\sm L=F_1\cup\cdots\cup F_n$ \st $\what\Upsilon\inv(\what M\sm
L)$ has at least $k$ \concompsns, and there exists an ends
decomposition $M\sm K=E_1\cup\cdots\cup E_m$ \st
$K\supset\Upsilon(L)$.
For each~$j=1,\dots,n$, part~(b) implies that $\Upsilon(F_j)\not\subset K$, and hence $F_j$
meets, and therefore contains, some \concomp of $\Upsilon\inv(M\sm
K)$. Thus, under the universal covering
$\Upsilon\circ\what\Upsilon\colon\wtil M\to M$, the inverse image
of $M\sm K$ has at least $k$ \concompsns, and therefore $\tilde
e(M)\geq k$. Thus $\tilde e(M)\geq\tilde e(\what M)$. Furthermore, if $\Upsilon$ is finite covering map, then the \concomps of the liftings of the ends in any ends
decompostion of $M$ form an ends decomposition for~$\what M$. Hence 
in this case we have $\tilde e(\what M)\geq\tilde e(M)$, and therefore
we have equality.
\end{pf}

\begin{defn}[cf.~Definition~2.2
of~\cite{NR Filtered ends}]\label{Special end definition} We will
call an end $E$ of a \con non\cpt complete Hermitian manifold $X$
\emph{special} if $E$ is of at least one of the following types:
\begin{enumerate}
\item[(BG)] $X$ has bounded geometry of order $2$ along $E$;

\item[(W)] There exists a \cont \plsh \fn $\vphi$ on~$X$ \st
\[
\setof{x\in E}{\vphi(x)<a}\Subset X\qquad\forall\,a\in\R;
\]

\item[(RH)] $E$ is a hyperbolic end and the Green's \fn vanishes
at infinity along $E$; or

\item[(SP)] $E$ is a parabolic end, the Ricci curvature of $g$ is
bounded below on $E$, and there exist positive constants $R$ and
$\delta$ \st \(\Vol\,\big(B(x;R)\big)>\delta\) for all $x\in E$.

\end{enumerate}
We will call an ends decomposition in which each of the ends is special a \emph{special ends decomposition}.
\end{defn}

According to \cite{Gro-Sur la groupe fond}, \cite{Li Structure
complete Kahler}, \cite{Gro-Kahler hyperbolicity},
\cite{Gromov-Schoen}, \cite{NR-Structure theorems},
\cite{Delzant-Gromov Cuts}, \cite{NR Filtered ends}, and \cite{NR L2 Castelnuovo}, a \con non\cpt complete K\"ahler manifold~$X$ that admits
a special ends decomposition and has at least three filtered ends admits a proper \holo mapping onto a
Riemann surface. One goal of this section is to show that if
$X$ has an irregular hyperbolic end of type~(BG), then two filtered ends suffice.
\begin{thm}\label{BG irregular hyp 2 ends thm}
If $X$ is a \con non\cpt hyperbolic complete K\"ahler manifold
that admits a special ends decomposition $X\sm K=E_1\cup\cdots\cup
E_m$ for which $E_1$ is an irregular hyperbolic end (i.e., $E_1$
contains an irregular sequence for $X$) of type~(BG) and $m\geq
2$, then $X$ admits a proper \holo mapping onto a Riemann surface.
\end{thm}
\begin{pf*}{Sketch of the proof}
Every end lying in a special end is itself special, so by the main
results of~\cite{NR-Structure theorems} and~\cite{NR Filtered
ends}, we may assume that $m=e(X)=2$. Moreover, as in the proof of
Theorem~3.4 of \cite{NR-Structure theorems}, we may also assume
that $E_2$ is a hyperbolic end of type~(BG). Theorem~2.6
of~\cite{NR-Structure theorems} then provides a nonconstant
bounded positive Dirichlet-finite \plh \fnns~$\rho_1$ on~$X$.
Proposition~\ref{plh fn infinite energy exist kahler prop} implies
that $X$ also admits a positive (quasi-Dirichlet-finite) \plh
\fnns~$\rho_2$ with bounded gradient and infinite energy. In
particular, the \holo $1$-forms $\theta_1\equiv\partial\rho_1$ and
$\theta_2\equiv\partial\rho_2$, are linearly independent, and
Theorems~0.1~and~0.2 of \cite{NR L2 Castelnuovo} give a proper
\holo mapping of~$X$ onto a Riemann surface.
\end{pf*}

\begin{lem}[cf.~Proposition~4.1 of~\cite{NR-BH Regular hyperbolic Kahler}]\label{cover to RS gives map to RS lem}
Let $(X,g)$ be a \con non\cpt complete K\"ahler manifold. If $X$
admits a special ends decomposition and some \con covering space
$\Upsilon\colon\what X\to X$ admits a proper \holo mapping onto a
Riemann surface, then $X$ admits a proper \holo mapping onto a
Riemann surface.
\end{lem}
\begin{pf}
The Cartan--Remmert reduction of $\what X$ is given by a proper
\holo mapping $\what\Phi\colon\what X\to\what S$ of $\what X$ onto
a Riemann surface~$\what S$ with $\what\Phi_*\ol_{\what
X}=\ol_{\what S}$. Fixing a fiber~$\what Z_0$ of~$\what\Phi$, we
may form a \rel \cpt \con \nbdns~$\what U_0$ of $\what Z_0$
in~$\what X$ and a nonnegative \cinf \plsh \fnns~$\hat\vphi_0$
on~$\what X\sm\what Z_0$ \st $\hat\vphi_0$ vanishes on $\what
X\sm\what U_0$ and $\vphi_0\to\infty$ at $\what Z_0$. The
image~$Z_0\equiv\Upsilon(\what Z_0)$ is then a \con \cpt \anal
subset of~$X$, and the \fn $\vphi_0\colon
x\mapsto\sum_{y\in\Upsilon\inv(x)}\hat\vphi_0(y)$ is a nonnegative
\cinf \plsh \fn on the domain $X\sm Z_0$ that vanishes on the
complement of the \rel \cpt \con \nbd~$U_0\equiv\Upsilon(\what
U_0)$ of~$Z_0$ in~$X$ and satisfies $\vphi_0\to\infty$ at~$Z_0$.

We may form a special ends decomposition $X\sm K=E_1\cup\cdots\cup
E_m$ with $K\supset U_0$, and setting $K_0\equiv K\sm
U_0$ and $E_0\equiv U_0\sm Z_0$, we get an ends decomposition
\[
(X\sm Z_0)\sm K_0=E_0\cup E_1\cup\cdots\cup E_m
\] 
of $X\sm Z_0$.
By part~(d) of Lemma~\ref{basic ends lem}, for $a\gg 0$, the set
$\setof{x\in X\sm Z_0}{\vphi(x)<a}$ has a \concompns~$Y_0$ that
contains $X\sm U_0$ and has the ends decomposition
\[
Y_0\sm K_0=E_0'\cup\cdots\cup E_m',
\]
where $E_j'\equiv E_j\cap Y_0$ for $j=0,\dots,m$. In particular,
the above is a special ends decomposition for the complete 
K\"ahler metric $g_0\equiv g+\lev{-\log(a-\vphi)}$ on~$Y_0$ (see,
for example, \cite{Dem}), with $E_0'$ regular hyperbolic and of
type~(W). Here, for any $\cal C^2$ \fnns~$\psi$, $\lev\psi$ denotes
the \emph{Levi form} of~$\psi$; that is, in local \holo coordinates $(z_1,\dots,z_n)$,
\[
\lev\psi=\sum_{j,k=1}^n\frac{\partial^2\vphi}{\partial z_j\partial\bar z_k}
dz_jd\bar z_k.
\]
Theorem~3.6 of~\cite{NR L2 Castelnuovo} implies that
there exists a nonconstant nonnegative \cont \plsh \fn on~$Y_0$
that vanishes on \(E_0'\cup K_0\) and, therefore, extends to a
\cont \plsh \fnns~$\alpha$ on~$X$ that vanishes on~$K$. 
Fixing a fiber $\what Z_1$ of $\what\Phi$ through a point at which
$\alpha\circ\Upsilon>0$, we see that, since $\alpha\circ\Upsilon$
is constant on~$\what Z_1$, the image $Z_1\equiv\Upsilon(\what
Z_1)$ must be a \con \cpt \anal subset of~$X\sm K\subset
Y_0\subset X\sm Z_0$. As above, we get a domain $Y_1\subset Y_0$,
a complete K\"ahler metric~$g_1$ on~$Y_1$, and a special ends
decomposition of $(Y_1,g_1)$ with at least three ends. Therefore,
by Theorem~3.4 of \cite{NR-Structure theorems} (or Theorem~3.1
of~\cite{NR Filtered ends}), there exists a proper \holo mapping
$\Phi_1\colon Y_1\to S_1$ of $Y_1$ onto a Riemann surface~$S_1$
\st $(\Phi_1)_*\ol_{Y_1}=\ol_{S_1}$. Forming the complement in~$X$
of two distinct fibers of~$\Phi_1$ and applying a construction
similar to the above, we get a proper \holo mapping $\Phi_2\colon
Y_2\to S_2$ of a domain $Y_2\supset X\sm Y_1$ in~$X$ onto a
Riemann surface~$S_2$ \st $(\Phi_2)_*\ol_{Y_2}=\ol_{S_2}$. The
maps $\Phi_1$ and $\Phi_2$ now determine a proper \holo
mapping~$\Phi$ of $X$ onto the Riemann surface
\[
S\equiv (S_1\sqcup
S_2)\big\slash\left[\Phi_1(x)\sim\Phi_2(x)\quad\forall\,x\in
Y_1\cap Y_2\right].
\]
\end{pf}

\begin{rmks}
1. The authors do not know whether or not the above lemma holds in general for the base an
arbitrary \con non\cpt complete K\"ahler manifold.

2. For the base a complete K\"ahler manifold with bounded geometry
(which is the relevant case for this paper), one may instead 
obtain the lemma from properness of the projection from the graph  over a suitable \irred \comp of the appropriate Barlet cycle space as 
in (Theorem~3.18 and the appendix of)~\cite{Campana}.
\end{rmks}

\begin{thm}\label{BG irreg hyp no BH cover thm}
Suppose $X$ is a \con non\cpt irregular hyperbolic complete
K\"ahler manifold with bounded geometry of order~$2$ and $e(X)=1$.
If $X$ admits a \con covering space $\Upsilon\colon\what X\to X$
with $H^1_c(\what X,\ol)\neq 0$, then $X$ admits a proper \holo
mapping onto a Riemann surface.
\end{thm}
\begin{pf}
Clearly, $\what X$ has bounded geometry of order~$2$. If $e(\what
X)\geq 3$ or $e(\what X)=1$, then $\what X$ admits a proper \holo
mapping onto a Riemann surface, and Lemma~\ref{cover to RS gives
map to RS lem} provides such a mapping on~$X$. Thus we may assume
that $e(\what X)=2$, and we may fix an ends decomposition $\what
X\sm K=E_1\cup E_2$. By part~(b) of Lemma~\ref{Behavior of
filtered ends for covering lem}, for $j=1,2$, $\Upsilon(E_j)$ is a
hyperbolic end of~$X$. It follows that $E_j$ is a hyperbolic end
of~$\what X$, since the lifting to~$\what X$ of a negative \cont
sub\harm \fn with supremum~$0$ on~$X$ is a negative \cont sub\harm
\fn on~$\what X$ with supremum~$0$ along~$E_j$.

Proposition~\ref{plh fn infinite energy exist kahler prop}
provides an unbounded positive \plh
\fnns~$\rho_1$ with bounded gradient and infinite energy on~$X$,
and we may set $\hat\rho_1\equiv\rho_1\circ\Upsilon$. Theorem~2.6
of~\cite{NR-Structure theorems} provides a nonconstant bounded \plh \fnns~$\hat\rho_2$ with finite energy on~$\what X$, and
Theorems~0.1 and 0.2 of \cite{NR L2 Castelnuovo}, applied to the
\holo $1$-forms $\partial\hat\rho_1$ and $\partial\hat\rho_2$,
give a proper \holo mapping of $\what X$ onto a Riemann surface.
Lemma~\ref{cover to RS gives map to RS lem} then gives the
required mapping on~$X$.
\end{pf}
Proposition~\ref{BH con boundary prop} provides some topological
conditions that give nonvanishing of the first compactly
supported cohomology with values in the structure sheaf.  In
particular, since any manifold with at least two filtered ends
admits a \con covering space with at least two ends, we get the
following consequence of Theorem~\ref{BG irreg hyp no BH cover thm}
(one may instead apply Theorem~\ref{BG irregular hyp 2 ends thm} 
and Lemma~\ref{cover to RS gives map to RS lem}):
\begin{cor}\label{BG irreg hyp 2 filtered end cor}
If $X$ is a \con non\cpt irregular hyperbolic complete K\"ahler
manifold with bounded geometry of order~$2$, $e(X)=1$, and $\tilde
e(X)\geq 2$, then $X$ admits a proper \holo mapping onto a Riemann
surface.
\end{cor}
We also get the following:
\begin{cor}\label{BG irreg hyp 2 pione boundary domain cor}
Suppose $X$ is a \con non\cpt irregular hyperbolic complete
K\"ahler manifold with bounded geometry of order~$2$, $e(X)=1$,
$\Omega$ is a nonempty smooth \rel \cpt domain in~$X$ for which
$E\equiv X\sm\overline\Omega$ is \con (i.e., $E$ is an end), and
$\Gamma'\equiv\text{\rm
im}\,[\pione{\overline\Omega}\to\pione{X}]$. If either
$\partial\Omega$ is not \conns, or $\partial\Omega$ is \con but
$\pione{\partial\Omega}$ does not surject onto $\Gamma'$, then $X$
admits a proper \holo mapping onto a Riemann surface.
\end{cor}
\begin{pf}
If $\partial\Omega$ is not \conns, then part~(a) of
Proposition~\ref{BH con boundary prop} implies that
$H^1_c(X,\ol)\neq 0$, and hence $X$ admits a proper \holo mapping
onto a Riemann surface. Suppose instead that
$C\equiv\partial\Omega$ is \conns, but
\(\Gamma\equiv\text{im}\,\left[\pione{C}\to\pione X\right]
\subsetneqq\Gamma'\). For a \con covering space
$\Upsilon\colon\what X\to X$ with $\Upsilon_*\pione{\what
X}=\Gamma$, $\Upsilon$ maps some \rel \cpt \con \nbdns~$U_0$ of
some \concompns~$C_0$ of $\what C\equiv\Upsilon\inv(C)$
isomorphically onto a \nbdns~$U$ of~$C$. By Theorem~\ref{BG irreg
hyp no BH cover thm}, we may assume that $e(\what X)=1$. The
unique \concompns~$\Omega_0$ of
$\what\Omega\equiv\Upsilon\inv(\Omega)$ for which $C_0$ is a
boundary \comp is a smooth domain, and
$C_0\subsetneqq\partial\Omega_0$. Moreover, each component of
$\what X\sm\overline\Omega_0$ must meet, and therefore contain, a
\comp of $\Upsilon\inv(E)$, so any such \comp must have non\cpt
closure. Proposition~\ref{BH con boundary prop} and
Theorem~\ref{BG irreg hyp no BH cover thm} together now give the
claim.
\end{pf}

\section{An irregular hyperbolic example}\label{BG irreg hyp example sect}

Because
the existence of irregular
hyperbolic complete K\"ahler manifolds with one end and bounded
geometry of order two is not completely obvious, 
an example is provided in this section. In fact, the
following is obtained:
\begin{thm}\label{Irregular hyp example thm}
There exists an irregular hyperbolic \con non\cpt complete
K\"ahler manifold~$X$ with bounded geometry of all orders \st
$e(X)=1$ and $\dim X=1$.
\end{thm}
\begin{rmk}
The authors do not know whether or not there exists an irregular
hyperbolic \con non\cpt complete K\"ahler manifold~$X$ with
bounded geometry of order~$0$ for which $H^1_c(X,\ol)=0$
(and hence which does not admit a proper \holo mapping onto a 
Riemann surface).
\end{rmk}

The idea of the construction is as follows. The complement of a
closed disk~$D$ in~$\C$ is irregular hyperbolic, but it has two
ends. Holomorphic attachment of a suitable sequence of tubes
(i.e., annuli) $\seq T\nu$, with boundary components $A_\nu$ and
$A'_\nu$ of $T_\nu$ for each~$\nu$, satisfying $A_\nu\to\infty$
and $A'_\nu\to p\in\partial D$, yields an irregular hyperbolic Riemann
surface with one end, and a direct construction yields a K\"ahler
metric with bounded geometry.
\begin{lem}\label{Irreg hyp ex disjoint disks lem}
Let $\set{\Delta(\zeta_\nu;R_\nu)}_{\nu=0}^\infty$ be a locally
finite sequence of disjoint disks in~$\C$. Then there exists a
sequence of positive numbers $\seq r\nu_{\nu=1}^\infty$ \st
$r_\nu<R_\nu$ for $\nu=1,2,3,\dots$, and
\[
b\equiv\sum_{\nu=1}^\infty\frac{\log\left[R_0\inv
R_\nu\inv(|\zeta_\nu-\zeta_0|+R_0)(|\zeta_\nu-\zeta_0|+R_\nu)\right]}
{\log\left[R_0\inv
r_\nu\inv(|\zeta_\nu-\zeta_0|+R_0)(|\zeta_\nu-\zeta_0|-r_\nu)\right]}<1.
\]
Moreover, for any such sequence $\seq r\nu$, the region
$\Omega\equiv\C\sm\bigcup_{\nu=1}^\infty
\overline{\Delta(\zeta_\nu;r_\nu)}$ is hyperbolic, and there
exists an irregular sequence $\seq\eta k$ in~$\Omega$ \st
$\eta_k\to\infty$ in~$\C$.
\end{lem}
\begin{pf}
It is easy to see that the above inequality will hold for all
sufficiently small positive sequences $\seq r\nu$. For each
$\nu=1,2,3,\dots$, let
\[
B_\nu\equiv\log\left[R_0\inv
R_\nu\inv(|\zeta_\nu-\zeta_0|+R_0)(|\zeta_\nu-\zeta_0|+R_\nu)\right],
\]
let
\[
C_\nu\equiv\log\left[R_0\inv
r_\nu\inv(|\zeta_\nu-\zeta_0|+R_0)(|\zeta_\nu-\zeta_0|-r_\nu)\right],
\]
and let $\alpha_\nu$ be the \harm \fn on
$\C\sm\set{\zeta_0,\zeta_\nu}$ given by
\[
z\mapsto\alpha_\nu(z)\equiv\frac
1{C_\nu}\log\left[\frac{|z-\zeta_0|}{|z-\zeta_\nu|}
\left(\frac{|\zeta_\nu-\zeta_0|+R_0}{R_0}\right)\right].
\]
Clearly, $B_\nu>0$, and since $|\zeta_\nu-\zeta_0|\geq R_\nu+R_0$,
$C_\nu>0$. At each point $z\in\partial\Delta(\zeta_0;R_0)$, we
have
\[
0\leq\alpha_\nu(z)=\frac 1{C_\nu}\log\left[
\frac{|\zeta_\nu-\zeta_0|+R_0}{|z-\zeta_\nu|}\right]\leq
\frac{B_\nu}{C_\nu},
\]
since $(|\zeta_\nu-\zeta_0|+R_\nu)|z-\zeta_\nu|\geq R_0R_\nu$. At
each point $z\in\partial\Delta(\zeta_\nu;r_\nu)$, we have
\[
\alpha_\nu(z)=\frac 1{C_\nu}\log\left[R_0\inv
r_\nu\inv(|\zeta_\nu-\zeta_0|+R_0)|z-\zeta_0|\right]\geq 1;
\]
while at each point $z\in\partial\Delta(\zeta_\nu;R_\nu)$, we have
\[
\alpha_\nu(z)=\frac 1{C_\nu}\log\left[R_0\inv
R_\nu\inv(|\zeta_\nu-\zeta_0|+R_0)|z-\zeta_0|\right]\leq
\frac{B_\nu}{C_\nu}.
\]
Moreover,
\[
0\leq\lim_{z\to\infty}\alpha_\nu(z)=\frac
1{C_\nu}\log\left[R_0\inv(|\zeta_\nu-\zeta_0|+R_0)\right]
\leq\frac{B_\nu}{C_\nu}.
\]
Therefore, since $\alpha_\nu$ is \harmns, we have $\alpha_\nu\geq
0$ on
\[
\C\sm\left[\Delta(\zeta_0;R_0)\cup \Delta(\zeta_\nu;r_\nu)\right]
\supset\overline\Omega\sm\Delta(\zeta_0;R_0),
\]
and $0\leq\alpha_\nu\leq\frac{B_\nu}{C_\nu}$ on
\(\C\sm\left[\Delta(\zeta_0;R_0)\cup\Delta(\zeta_\nu;R_\nu)\right]\).
Consequently, the series $\sum_{\nu=1}^\infty\alpha_\nu$ converges
uniformly on \cpt subsets
of~$\overline\Omega\sm\Delta(\zeta_0;R_0)$ to a nonnegative \cont
\fnns~$\alpha$ \st $\alpha$ is positive and \harm
on~$\Omega\sm\overline{\Delta(\zeta_0;R_0)}$, $\alpha\leq b<1$ on
the set
\[
I\equiv\Omega\sm\bigcup_{\nu=0}^\infty\Delta(\zeta_\nu;R_\nu)
=\C\sm\bigcup_{\nu=0}^\infty \Delta(\zeta_\nu;R_\nu),
\]
and for each~$\nu=1,2,3,\dots$, we have $0<\alpha-\alpha_\nu<1$ on
$\overline{\Delta(\zeta_\nu;R_\nu)}\sm\Delta(\zeta_\nu;r_\nu)$ and
$\alpha>\alpha_\nu\geq 1$ on $\partial\Delta(z_\nu;r_\nu)$.

Clearly, $\Omega\supset\overline{\Delta(\zeta_0;R_0)}$ is
hyperbolic, and the \hmib $\Omega$ with respect to
$\Omega\sm\overline{\Delta(\zeta_0;R_0)}$ extends to a \cont \fn
$u\colon\overline\Omega\sm\Delta(\zeta_0;R_0)\to[0,1]$. For each
$R>R_0$, the \cont \fn $\beta_R$ on
$\overline\Omega\sm\Delta(\zeta_0;R_0)$ given by
\[
z\mapsto\beta_R(z)\equiv\alpha(z)
+\frac{\log\left(|z-\zeta_0|/R_0\right)}{\log(R/R_0)}
\]
is \harm on $\Omega\sm\overline{\Delta(\zeta_0;R_0)}$ and
satisfies $\beta_R>\alpha>1=u$ on
$\bigcup_{\nu=1}^\infty\partial\Delta(\zeta_\nu;r_\nu)$,
$\beta_R=\alpha\geq 0=u$ on $\partial\Delta(\zeta_0;R_0)$, and
$\beta_R\geq 1\geq u$ on
$\overline\Omega\cap\partial\Delta(\zeta_0;R)$. Hence $\beta_R\geq
u$ on
$(\overline\Omega\sm\Delta(\zeta_0;R_0))\cap\overline{\Delta(\zeta_0;R)}$.
Passing to the pointwise limit as $R\to\infty$, we get $\alpha\geq
u$ on $\overline\Omega\sm\Delta(\zeta_0;R_0)$. However,
$\alpha\leq b<1$ on the set
\(I\subset\Omega\sm\Delta(\zeta_0;R_0)\),
so any sequence $\seq\eta k$ in~$I$ with $\eta_k\to\infty$ in~$\C$
is an irregular sequence in~$\Omega$.
\end{pf}

\begin{lem}\label{BG from fn on C lem}
Let $k$ be a positive integer, and $\rho$ a positive \cinf \fn
on~$\C$ \st $D\rho,D^2\rho,D^3\rho,\dots,D^k\rho$ are bounded.
Then the complete K\"ahler metric $g\equiv e^{2\rho}g_{\C}$ has
bounded geometry of order~$k$. In fact, the pullbacks of~$g$ under
the local \holo charts
\[
\Psi_{z_0}\colon\Delta(0;1)\to\Delta(z_0;e^{-\rho(z_0)})
\]
given by $\Psi_{z_0}\colon z\mapsto e^{-\rho(z_0)}z+z_0$, for each
point $z_0\in\C$, have the appropriate uniformly bounded
derivatives.
\end{lem}
\begin{pf}
For each point $z_0\in\C$, the pullback of the associated
$(1,1)$-form $\omega_g\equiv e^{2\rho}\frac i2dz\wedge d\bar z$
under~$\Psi_{z_0}$ is given by
\[
\Psi_{z_0}^*\omega_g
=e^{2\left(\rho(\Psi_{z_0})-\rho(z_0)\right)}\frac i2dz\wedge
d\bar z.
\]
The bound on $D\rho$ gives a Lipschitz constant~$C$ for~$\rho$,
and hence
\[
e^{-2C}\leq e^{-2Ce^{-\rho(z_0)}}\leq
e^{2\left(\rho(\Psi_{z_0})-\rho(z_0)\right)}\leq
e^{2Ce^{-\rho(z_0)}}\leq e^{2C}.
\]
A similar argument gives uniform bounds on the $m$th order
derivatives of the \fns
$\left\{e^{2\left(\rho(\Psi_{z_0})-\rho(z_0)\right)}\right\}_{z_0\in\C}$
for $m=1,\dots,k$.
\end{pf}

\begin{pf*}{Proof of Theorem~\ref{Irregular hyp example thm}}
\textbf{Step 1.} Construction of a suitable irregular hyperbolic
region in~$\C$. Let us fix a constant~$R>1$ and disjoint disks
$\set{\Delta(\zeta_\nu;R)}_{\nu=0}^\infty$ \st $\zeta_0=0$ and
$\zeta_\nu\to\infty$ so fast that
\[
\sum_{\nu=1}^\infty\frac{\log\left[R^{-2}(|\zeta_\nu|+R)^2\right]}
{\log\left[R\inv
e^{|\zeta_\nu|}(|\zeta_\nu|+R)(|\zeta_\nu|-e^{-|\zeta_\nu|})\right]}<1.
\]
In particular, by Lemma~\ref{Irreg hyp ex disjoint disks lem}, the
domain
\[
\Omega_0\equiv\C\sm\bigcup_{\nu=1}^\infty
\overline{\Delta(\zeta_\nu;e^{-|\zeta_\nu|})}
\]
is irregular hyperbolic; in fact, there exists an irregular
sequence $\seq\eta k$ in~$\Omega_0$ \st $\eta_k\to\infty$ in~$\C$.

\textbf{Step 2. Construction of a bounded geometry K\"ahler metric
on a region.} By Lemma~\ref{BG from fn on C lem}, fixing a
positive \cinf \fnns~$\rho$ on~$\C$ \st $\rho(z)=|z|$ on a \nbd of
$\C\sm\Delta(0;1)$, we get a complete K\"ahler metric $g_0\equiv
e^{2\rho}g_{\C}$ with bounded geometry of all orders on~$\C$ and
associated local \holo charts
\[
\Psi_{z_0}\colon\Delta(0;1)\to\Delta(z_0;e^{-\rho(z_0)})
\]
given by $\Psi_{z_0}\colon z\mapsto e^{-\rho(z_0)}z+z_0$ for each
point $z_0\in\C$. Letting $R_0$ and $R_1$ be constants with
$1<R_0<R_1<R$, $g_{\mathbb H}$ the standard hyperbolic metric on
the upper half plane~$\mathbb H$, $\Phi$ a M\"obius transformation
with $\Phi((\C\sm\overline{\Delta(0;1)})\cup\set\infty)=\mathbb H$
and $\imag\Phi>5R$ on $(\C\sm\Delta(0;R_0))\cup\set\infty$,
$g_1\equiv\Phi^*g_{\mathbb H}$, and $\lambda\colon\C\to[0,1]$ a
\cinf \fn with $\lambda\equiv 0$ on  a \nbd of
$\overline{\Delta(0;R_0)}$ and $\lambda\equiv 1$ on
$\C\sm\Delta(0;R_1)$, we get a complete K\"ahler metric
\[
g_2\equiv\lambda g_0+(1-\lambda)g_1
\]
with bounded geometry of all orders on the region
\(\C\sm\overline{\Delta(0;1)}\). Setting $\xi_\nu\equiv 2\nu
R+i2R$ for each~$\nu=1,2,3,\dots$, we get disjoint disks
$\set{\Delta(\xi_\nu;R)}_{\nu=1}^\infty$ in $\setof{z\in\mathbb
H}{R<\imag z<3R}$ (and an isometric isomorphism
$\Delta(\xi_1;R)\to\Delta(\xi_\nu;R)$ in $\mathbb H$ given by
$z\mapsto z+2(\nu-1)R$ for each~$\nu$). We have
\[
\overline{\Delta(0;1)}
\cup\bigcup_{\nu=1}^\infty\Phi\inv(\overline{\Delta(\xi_\nu;R)})
\subset\Delta(0;R_0)\Subset\Delta(0;R)\Subset\Omega_0,
\]
and hence we have a region
\[
\Omega_1\equiv\Omega_0\sm\left(\overline{\Delta(0;1)}
\cup\bigcup_{\nu=1}^\infty\Phi\inv(\overline{\Delta(\xi_\nu;1)})\right).
\]

\textbf{Step 3. Construction of the Riemann surface $X$.} For each
$\nu=1,2,3,\dots$, let $T_\nu$ be a copy of the annulus
\(\Delta(0;1/R,R)\equiv \setof{z\in\C}{1/R<|z|<R}\), and let
\[
\Lambda_\nu\colon\C\to\C \qquad\text{and}\qquad
\Upsilon_\nu\colon\C^*\to\C^*
\]
be the biholomorphisms given by $w\mapsto
e^{-|\zeta_\nu|}w+\zeta_\nu$ and $w\mapsto\frac 1w+\xi_\nu$,
respectively. We then get a Riemann surface
\[
X\equiv\left(\Omega_1\sqcup\bigsqcup_{\nu=1}^\infty
T_\nu\right)\bigg/{\sim},
\]
where for each $\nu=1,2,3,\dots$, and each $w\in T_\nu$,
$z\in\Delta(\zeta_\nu;e^{-|\zeta_\nu|},Re^{-|\zeta_\nu|})$
satisfies
\[
z\sim w\qquad\iff\qquad z=\Lambda_\nu(w),
\]
and $z\in\Phi\inv(\Delta(\xi_\nu;1,R))$ satisfies
\[
z\sim w\qquad\iff\qquad \Phi(z)=\Upsilon_\nu(w).
\]
$X$ is hyperbolic, because for each point
$z_0\in(\partial\Delta(0;1))\sm\set{\Phi\inv(\infty)}
\subset\partial\Omega_1$,
there exists a barrier~$\beta$ on~$\Omega_1$ at~$z_0$ and a \rel
\cpt \nbdns~$U$ of $z_0$ in~$\C$ \st $\overline
U\sm\overline{\Delta(0;1)}\subset\Omega_1$ and $\beta$ is equal
to~$-1$ on $\Omega_1\sm U$, and thus we may extend~$\beta$ to a
\cont sub\harm \fn on~$X$ that is equal to~$-1$ on
$X\sm(\Omega_1\cap U)$.  Fixing a disk
\[
D\Subset\Delta(0;R)\cap\Omega_1\subset X,
\]
and letting $u\colon X\sm D\to[0,1)$ be the \hmib $X$ with respect
to $X\sm\overline D$, we see that the restriction
$u\restrict{\Omega_0\sm\Delta(0;R)}$ cannot approach~$1$ along the
sequence $\seq\eta k$, so $X$ must be irregular hyperbolic. It is
easy to see that $e(X)=1$.

\textbf{Step 4. Construction of a bounded geometry K\"ahler metric
on~$X$.} Let us fix a \cinf \fnns~$\tau$ on~$\C$ \st
$0\leq\tau\leq 1$, $\tau\equiv 1$ on $\C\sm\Delta(0;R_0)$, and
$\tau\equiv 0$ on $\Delta(0;1/R_0)$. Then we get a K\"ahler
metric~$g$ on~$X$ by setting $g=g_2$ on
\[
\Omega_2\equiv\C\sm\left(\overline{\Delta(0;1)}\cup
\bigcup_{\nu=1}^\infty
\overline{\Delta(\zeta_\nu;R_0e^{-|\zeta_\nu|})} \cup
\bigcup_{\nu=1}^\infty\Phi\inv(
\overline{\Delta(\xi_\nu;R_0)})\right)\subset\Omega_1\subset X,
\]
and \(g=\tau\Lambda_\nu^*g_0+(1-\tau)\Upsilon_\nu^*g_{\mathbb H}\)
on~$T_\nu\subset X$ for each $\nu=1,2,3,\dots$.

For $\nu=1,2,3,\dots$, on $T_\nu$ we have
\(\Lambda_\nu^*g_0=e^{2(|e^{-|\zeta_\nu|}w+\zeta_\nu|-|\zeta_\nu|)}g_{\C}\)
and \(\Upsilon_\nu^*g_{\mathbb H}=\Upsilon_1^*g_{\mathbb H}\)
(since \(\Upsilon_\nu=\Upsilon_1+2(\nu-1)R\)).
Therefore, since the \fns 
\[
w\mapsto
|e^{-|\zeta_\nu|}w+\zeta_\nu|-|\zeta_\nu|\in[-Re^{-R},Re^{-R}]\qquad\text{for }\nu=1,2,3,\dots,
\]
have uniformly bounded derivatives of order~$k$
on~$T_\nu=\Delta(0;1/R,R)$ for each~$k=0,1,2,\dots$, $(X,g)$ has
bounded geometry of all orders along $X\sm\Omega_3$, where
\[
\Omega_3\equiv\C\sm\left(\overline{\Delta(0;1)}\cup
\bigcup_{\nu=1}^\infty
\overline{\Delta(\zeta_\nu;R_1e^{-|\zeta_\nu|})} \cup
\bigcup_{\nu=1}^\infty\Phi\inv(
\overline{\Delta(\xi_\nu;R_1)})\right)\subset\Omega_2\subset\Omega_1\subset
X.
\]

There exists a positive constant $r_0$ \st for each point
$z_0\in\Omega_3\cap\Delta(0;R_0)$, we have
\[
B\equiv B_{g_{\mathbb
H}}(\Phi(z_0);r_0)\subset\Phi(\Omega_2\cap\Delta(0;R_1))
\]
and $g=g_1=\Phi^*g_{\mathbb H}$ on
$\Phi\inv(B)\subset\Omega_2\cap\Delta(0;R_1)$. Thus $(X,g)$ has
bounded geometry of all orders 
along $\Omega_3\cap\Delta(0;R_0)$, as well as
along the \cpt set $\overline{\Delta(0;R_0,R)}\subset\Omega_1$.

Finally, if $r_1$ is a constant with \(0<r_1<\min(1,R-R_1)\), and
$z_0\in\Omega_3\sm\Delta(0;R)$, then
\[
\Delta(z_0;r_1e^{-\rho(z_0)})\cap\Delta(0;R_1)=\emptyset.
\]
Moreover, if $z\in\Delta(z_0;r_1e^{-\rho(z_0)})\cap
\overline{\Delta(\zeta_\nu;R_0e^{-|\zeta_\nu|})}$ for some~$\nu$, then
\begin{align*}
R_1e^{-|\zeta_\nu|}&<|z_0-\zeta_\nu|<r_1e^{-|z_0|}+R_0e^{-|\zeta_\nu|}
\leq(r_1e^{|\zeta_\nu-z_0|}+R_0)e^{-|\zeta_\nu|}\\
&\leq(r_1e^{r_1e^{-|z_0|}+R_0e^{-|\zeta_\nu|}}+R_0)e^{-|\zeta_\nu|}
\leq(r_1e^{r_1e^{-R}+R_0e^{-R}}+R_0)e^{-|\zeta_\nu|}.
\end{align*}
Thus for $r_1$ sufficiently small, we will have, for every point
$z_0\in\Omega_3\sm\Delta(0;R)$,
\[
D_{z_0}\equiv\Delta(z_0;r_1e^{-\rho(z_0)})\subset\Omega_2\sm\Delta(0;R_1),
\]
and in particular, $g=g_2=g_0$ on~$D_{z_0}$. The resulting family
of biholomorphisms  $\Delta(0;1)\to D_{z_0}$ given by $z\mapsto
r_1ze^{-|z_0|}+z_0$ for each such point~$z_0$ then have the
required uniform bounds, so $(X,g)$ has bounded geometry of all
orders along $\Omega_3\sm\Delta(0;R)$, and therefore along~$X$
itself, and completeness follows.
\end{pf*}

\bibliographystyle{amsalpha.bst}

\end{document}